\documentclass[a4paper,12pt,reqno]{amsart}
\usepackage{amsmath,graphicx,graphics,amssymb}
\newtheorem{lem}{Lemma}

\newtheorem{theo}{Theorem}

\newtheorem{cor}{Corollary}

\newtheorem{prop}{Proposition}

\numberwithin{equation}{section}

\newcommand{\M}{\operatorname{M}}
\newcommand{\W}{\operatorname{W}}

\newcommand{\epf}{\hfill{$\square$}\medskip}
\newcommand{\Z}{\mathbb Z}

\newmuskip\pFqskip
\mathchardef\pFcomma=\mathcode`, 

\begin{document}

\title{Symmetries of shamrocks IV: The self-complementary case}

\author{Mihai Ciucu}
\address{Department of Mathematics, Indiana University, Bloomington, Indiana 47405}

\thanks{Research supported in part by NSF grant DMS-1501052}

\begin{abstract} In this paper we enumerate the centrally symmetric lozenge tilings of a hexagon with a shamrock removed from its center. Our proof is based on a variant of Kuo's graphical condensation method in which only three of the four involved vertices are on the same face. As a special case, we obtain a new proof of the enumeration of the self-complementary plane partitions.
\end{abstract}

\maketitle

\section{Introduction}



This paper continues our program of enumerating the symmetry classes of lozenge tilings of the regions we introduced in \cite{vf}, where we were motivated by generalizing MacMahon's beautiful formula \cite{MacM} stating that the number of lozenge tilings\footnote{ A lozenge tiling of a lattice region $R$ on the triangular lattice is a covering of $R$ by lozenges (i.e. unions of unit triangles that share an edge), with no gaps or overlaps.} of a hexagonal region of side-lengths $x$, $y$, $z$, $x$, $y$, $z$ (in cyclic order) on a triangular lattice is
\begin{equation}
P(x,y,z)=\prod_{i=1}^x\prod_{j=1}^y\prod_{k=1}^z \frac{i+j+k-1}{i+j+k-2}.
\label{eaa}
\end{equation}
(MacMahon phrased this result in terms of plane partitions --- gravitationally stable stacks of unit cubes in a corner --- that fit in a given box; see e.g. \cite{DT} for the equivalence of these to lozenge tilings).

In \cite{vf} we showed that the number of lozenge tilings of hexagonal regions with a 4-lobed structure (called a shamrock) removed from their center is given by a product formula generalizing \eqref{eaa} (see \cite{fv} for a counterpart of this). Motivated by the truly remarkable fact that all symmetry classes of lozenge tilings of a hexagon are given by equally beautiful formulas (see \cite{And}\cite{Sta}\cite{Kup}\cite{Ste}\cite{KKZ} and the survey \cite{Bres} for more recent developments), we consider the problem of enumerating the symmetry classes of tilings of these more general regions. Six new questions arise in this way. In \cite{symffa} we solved the cyclically symmetric and the cyclically symmetric and transpose complementary cases (invariance under rotation by 120 degrees, resp. invariance under the same plus reflection across vertical), in \cite{symffb} the transpose complementary case (invariance under reflection across the vertical), and in  \cite{symffc} the symmetric and self-complementary case (tilings that are both vertically symmetric and centrally symmetric).

The purpose of this paper is to present the enumeration of a fifth case, that of self-com\-ple\-men\-ta\-ry tilings (i.e., tilings that are centrally symmetric). In order for such tilings to exist, the region must itself-be centrally symmetric. This implies that two of the lobes of the shamrock are empty, and the third has the same size as the core. This leaves the case of horizontally symmetric tilings of our regions, which will be presented in a separate paper.

Our results (see Theorems \ref{tbb} and \ref{tbc}) --- for which we provide simple, combinatorial proofs --- can also be viewed as providing a generalization of the enumeration of self-complementary plane partitions that fit in a box, which was first proved by Stanley \cite{Sta}. 

It turns out that a key role for proving such formulas is played by Kuo's graphical condensation method \cite{KuoOne}\cite{KuoTwo} (see also \cite{genKuo} for a generalization, and \cite{CFone}\cite{CFtwo}\cite{CFthree} for more applications). However, for the regions involved the symmetry class under consideration we run into the following difficulty: Kuo condensation works best when the involved four vertices can be chosen so that, upon removing an even number of them, one is led to graphs of the same type; it turns out that for the case at hand this basically determines the positions of these four vertices; the difficulty we face is that the resulting four vertices are not on the same face of the graph.

Generally speaking, this should preclude us from using graphical condensation altogether, as all the available results in the literature assume the involved vertices to be on the same face. The striking fact is that for a large family of centrally symmetric graphs,
a related, unexpected recurrence turns out to hold (see Theorems \ref{tba}, \ref{tbaa} and \ref{tbab}), and this allows us to prove our results.

\section{Statement of main results}

For a graph $G$ with weights on its edges, we denote by $\M(G)$ the sum of the weights\footnote{ The weight of a matching is the product of the weights of the edges in it.} of all the perfect matchings of $G$. Clearly, setting all the edge weights equal to 1 results in $\M(G)$ simply counting all such perfect matchings.

Before stating the formulas that enumerate centrally symmetric tilings of our regions, we present a surprising and unusual Kuo-style recurrence we found, which allows us to prove them.

\begin{theo}
\label{tba}
Let $G$ be a weighted $($not necessarily planar!$\,)$ graph. Let $a$, $b$, $c$ and $d$ be four arbitrary, distinct vertices of $G$. Assume that for all two-element subsets $S\subset\{a,b,c,d\}$, all paths connecting $a$ to $b$, $a$ to $c$, or $b$ to $c$ that arise from the superposition of any perfect matching $\mu$ of $G\setminus S$ with any perfect matching  $\nu$ of $G\setminus\overline{S}$ have  odd length\footnote{ Here $\overline{S}$ denotes the complement of $S$ in the set $\{a,b,c,d\}$.}.

Then we have
\begin{align}
&\!\!\!\!\!\!
\M(G)\M(G\setminus\{a,b,c,d\})
=
\M(G\setminus\{a,b\})\M(G\setminus\{c,d\})
+
\M(G\setminus\{a,c\})\M(G\setminus\{b,d\})
\nonumber
\\
&\ \ \ \ \ \ \ \ \ \ \ \ \ \ \ \ \ \ \ \ \ \ \ \ \ \ \ \ \ \ \ \ \ 
+
\M(G\setminus\{a,d\})\M(G\setminus\{b,c\}).
\label{eba}
\end{align}
\end{theo}


\parindent0pt
\textsc{Remark 1.} Both the hypothesis and the statement \eqref{eba} of the above theorem are visibly symmetric in $a,b,c,d$, with the exception of the assumption on the oddness of the lengths of the described paths, as their endpoints must be in the set $\{a,b,c\}$. However, it is not hard to see that as long as $G$ has an even number of vertices
, the stated $\{a,b,c\}$-version of the condition is equivalent to each of its $\{a,b,d\}$-, $\{a,c,d\}$- and $\{b,c,d\}$-versions. 
\parindent15pt  

\medskip
\parindent0pt
\textsc{Remark 2.} The above recurrence has a striking resemblance to Kuo's original graphical condensation formula \cite{KuoTwo} for non-bipartite graphs. Indeed, the latter --- which recall holds when $a$, $b$, $c$, $d$ occur in this cyclical order on the same face of $G$ --- is obtained from \eqref{eba} by simply changing the sign of the middle term on the right hand side!
\parindent15pt  

\medskip
We present the proof of Theorem \ref{tba} in the next section. Simple combinatorial arguments which we present in Section 4 will imply that the condition on the oddness of the lengths of paths in the hypothesis of Theorem \ref{tba} holds in the natural context described in Theorem \ref{tbaa} below, and thus the latter will follow from Theorem \ref{tba}.

\begin{theo}
\label{tbaa}
Let $G$ be a weighted, centrally symmetric, planar bipartite graph embedded in an annulus. Let $\{a_1,a_2\}$,  $\{b_1,b_2\}$ and $\{c_1,c_2\}$ be pairs of symmetric vertices on the outer face of~$G$. Assume that $a_1$, $b_1$, $c_1$, $a_2$, $b_2$, $c_2$ appear in this cyclic order around the outer face, and that they alternate in color. Let $\{d_1,d_2\}$ be a pair of symmetric vertices of $G$ on the central face. Then we have
\begin{align}
\M_{\odot}(G)\M_{\odot}(G_{abcd})=\M_{\odot}(G_{ab})\M_{\odot}(G_{cd})+
\M_{\odot}(G_{ac})\M_{\odot}(G_{bd})+
\M_{\odot}(G_{ad})\M_{\odot}(G_{bc}),
\label{ebaa}
\end{align}
where $G_{abcd}=G\setminus\{a_1,a_2,b_1,b_2,c_1,c_2,d_1,d_2\}$, $G_{ab}=G\setminus\{a_1,a_2,b_1,b_2\}$, etc., and $\M_{\odot}(H)$ denotes the number of centrally symmetric perfect matchings of the graph $H$. 

\end{theo}

We will also need a slight extension of the above result, which we present in Theorem \ref{tbab} below. Suppose that $G$ is embedded in a disk, and allow the symmetric vertices $d_1$ and $d_2$ to be on two different symmetric faces $F_1$ and $F_2$. To state our extension of Theorem \ref{tbaa} we need the following definition.

Denote by $\mathcal{M}_{\odot}(H)$ the set of centrally symmetric perfect matchings of the graph $H$. Let $x\in\{a,b,c\}$, and consider the centrally symmetric matchings $\mu\in\mathcal{M}_{\odot}(G_{xd})$ and $\nu\in\mathcal{M}_{\odot}(G_{\overline{xd}})$.
In the superposition of $\mu$ and $\nu$ there is a unique pair of symmetric paths $(P_1,P_2)$ so that each $P_i$ connects vertices from among $\{a_1,a_2,b_1,b_2,c_1,c_2\}$. We say that $(\mu,\nu)$ is {\it special} if each of these $P_i$'s has endpoints of the same color, and passes in between the faces $F_1$ and $F_2$.

\begin{theo}
\label{tbab}
Let $G$ be a centrally symmetric, weighted planar bipartite graph embedded in a disk. Let $\{a_1,a_2\}$,  $\{b_1,b_2\}$ and $\{c_1,c_2\}$ be pairs of symmetric vertices on the outer face of $G$. Assume that $a_1$, $b_1$, $c_1$, $a_2$, $b_2$, $c_2$ appear in this cyclic order around the outer face, and that they alternate in color. Let $\{d_1,d_2\}$ be a pair of symmetric vertices of $G$ on two symmetric faces $F_1$ and $F_2$ .
Then we have
\begin{align}
&
  \M_{\odot}(G)\M_{\odot}(G_{abcd})=\M_{\odot}(G_{ab})\M_{\odot}(G_{cd})+
\M_{\odot}(G_{ac})\M_{\odot}(G_{bd})+
\M_{\odot}(G_{ad})\M_{\odot}(G_{bc})
\nonumber
\\
&\ \ \ \ \ \ \ \ \ \ \ \ \ \ \ \ 
-\W\left(\bigcup_{x\in\{a,b,c\}}\left\{(\mu,\nu):\mu\in\mathcal{M}_{\odot}(G_{xd}),
\nu\in\mathcal{M}_{\odot}(G_{\overline{xd}}),(\mu,\nu)\ {\operatorname{special}}\right\}\right),
\label{ebab}
\end{align}
where for a set $A$ of pairs of perfect matchings, $\W(A)$ is the sum of the weights of the elements of $A$. 
\end{theo}

\parindent0pt
\textsc{Remark 3}. The above extensions of Kuo condensation to centrally symmetric graphs bear some resemblance to Kuperberg's extension (see  Theorem 6 in \cite{Kupex}) of the Kasteleyn-Temperley-Fisher theorem (see \cite{Kast}\cite{TF}) on the enumeration of perfect matchings of a planar graph by a Pfaffian to graphs embedded in the projective plane.

\parindent15pt

We now turn to stating our formulas that enumerate centrally symmetric tilings of our regions.

For non-negative integers $x$, $y$, $z$ of the same parity, let $B_{x,y,z,k}$ be the region obtained from the hexagon $H_{x,y,z}$ of side-lengths $x$, $y$, $z$, $x$, $y$, $z$ (in cyclic order) on the triangular lattice, by removing from its center a bowtie with lobes of size $k$; we assume (as this can be done without loss of generality) that the edges of the bowtie that do not touch the center are parallel to the sides of length $z$. It is not hard to see that in order for $B_{x,y,z,k}$ to have centrally symmetric tilings, $k$ must have the same parity as $x$, $y$ and $z$.

Strictly speaking, handling the regions $B_{x,y,z,k}$ covers the self-complementary symmetry class of hexagons with a shamrock removed from the center, which was the original motivation for the current paper. However, from the point of view of obtaining (by specializing $k=0$) a new proof of the enumeration of self-complementary plane partitions, this only covers the case when all three pairs of opposite edges of the hexagon have even side-lengths, and not the cases when two of these are are even and one odd, or when two are odd and one even (there are no centrally symmetric tilings of a hexagon with all three side-lengths odd).

To complete this goal, we introduce another family of regions as follows. Let $x$, $y$ and $z$ be  non-negative integers so that $x$ and $z$ have the same parity, while $y$ has the opposite parity. Then the center of the hexagon $H_{x,y,z}$ is the midpoint of a lattice segment $s$ parallel to the polar direction $2\pi/3$. Define $B'_{x,y,z,k}$ to be the region obtained by removing from  $H_{x,y,z}$ a ``disconnected bowtie'' consisting of two equilateral triangles of side-length $k$ placed so that the segment $s$ joins their closest vertices. Unlike the regions $B_{x,y,z,k}$, it turns out that the $B'_{x,y,z,k}$'s have centrally symmetric tilings for both parities of $k$.

It is not hard to see that unless $k\leq x,y,z$, neither $B_{x,y,z,k}$ nor $B'_{x,y,z,k}$ have any lozenge tilings. We can therefore assume without loss of generality that $k\leq x,y,z$ in our formulas for the number of lozenge tilings of the regions $B_{x,y,z,k}$ and $B'_{x,y,z,k}$.

\begin{figure}[h]
  \centerline{
\hfill
{\includegraphics[width=0.44\textwidth]{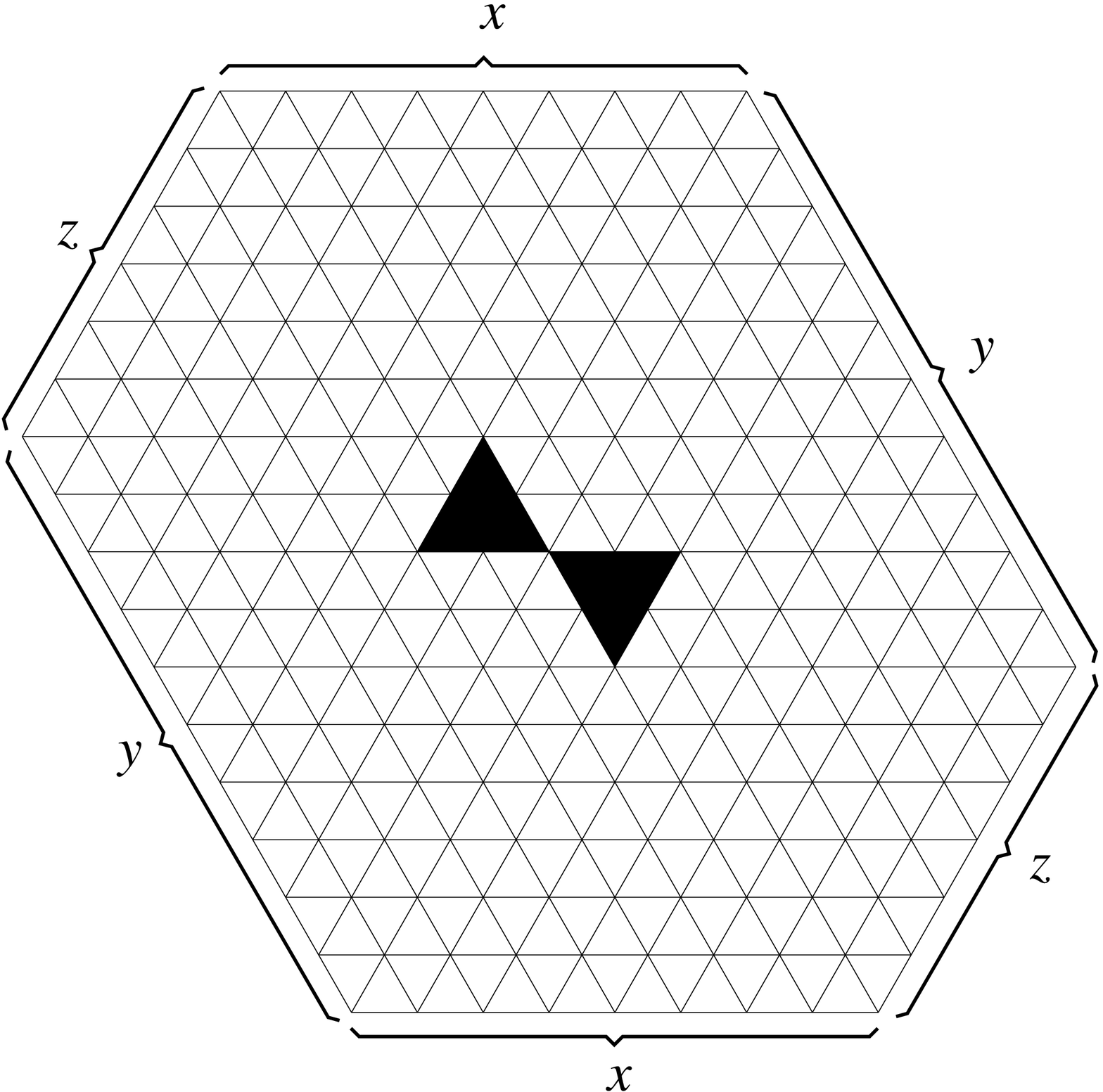}}
\hfill
{\includegraphics[width=0.44\textwidth]{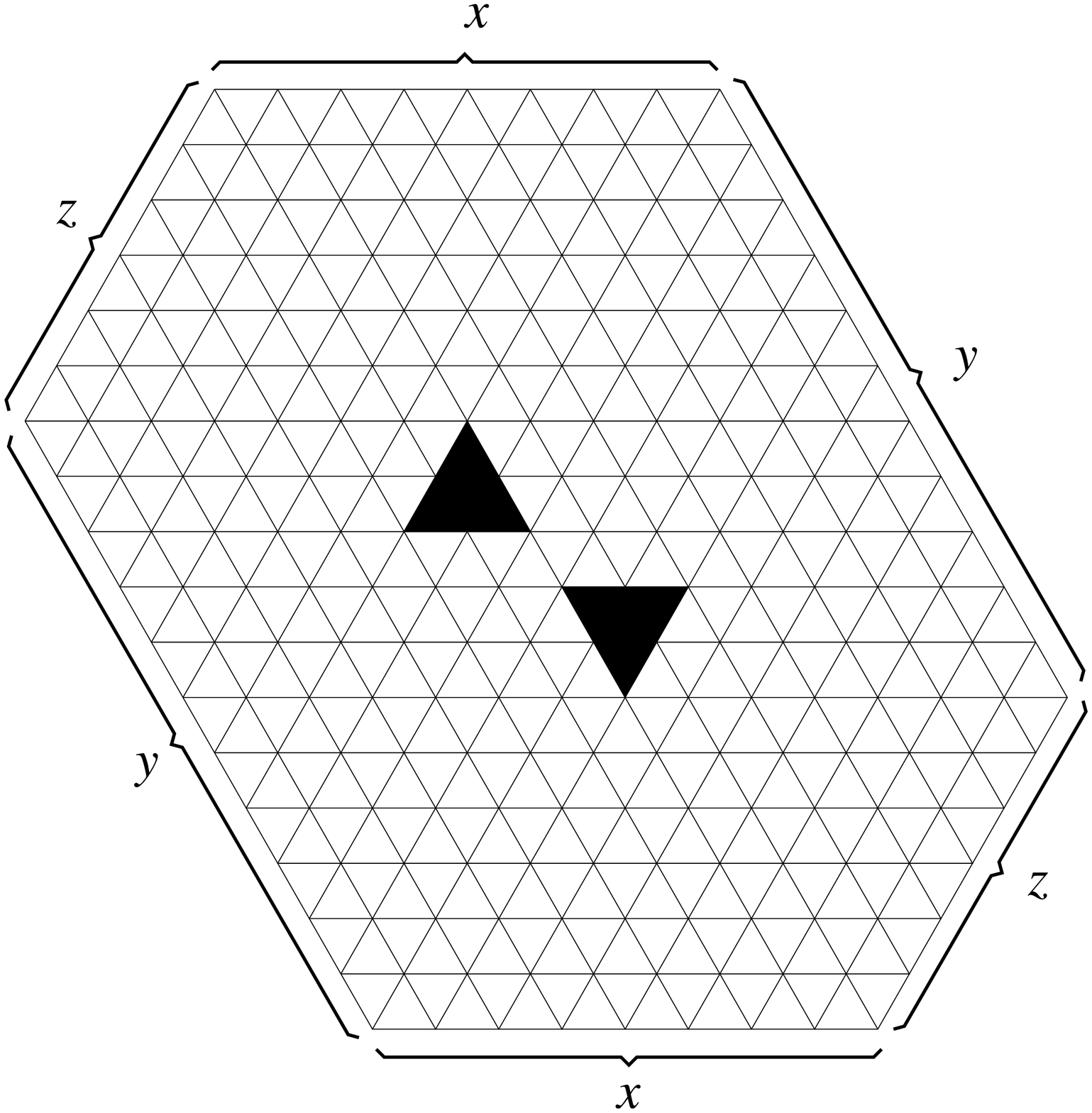}}
\hfill
}
\caption{\label{fba} The regions $B_{x,y,z,k}$ for $x=8, y=10, z=6, k=2$ (left) and $B'_{x,y,z,k}$ for $x=8, y=11, z=6, k=2$ (right).}
\label{fbb}
\end{figure}

The main results of this paper are the following. Recall that $P(x,y,z)$ is given by \eqref{eaa}.

\begin{theo}
\label{tbb}
For non-negative integers $x$, $y$, $z$ and $k$ of the same parity, with $k\leq x,y,z$, we have
\begin{align}
\M_{\odot}(B_{x,y,z,k})=P\left(\frac{y+k}{2},\frac{z-k}{2},k\right)
\prod_{i=1}^{(y-k)/2}\frac{\left(\frac{x-k}{2}+i\right)_k}{(i)_k}
\prod_{i=1}^{(z-k)/2}\frac{\left(\frac{x+k}{2}+i\right)_{(y-k)/2}}{(k+i)_{(y-k)/2}}
\frac{\left(\frac{x+k}{2}+i\right)_{(y+k)/2}}{(k+i)_{(y+k)/2}}.
\label{ebb}
\end{align}
\end{theo}

\begin{theo}
\label{tbc}
$(${\rm a}$)$. Let $x$, $y$, $z$, $k$ be non-negative integers, $k\leq x,y,z$, and assume that $x$, $z$, and $k$ have parity opposite to the parity of $y$. Then
\begin{align}
&
\M_{\odot}(B'_{x,y,z,k})=P\left(\frac{y+k+1}{2},\frac{z-k}{2},k\right)
\nonumber
\\
&\ \ \ \ \ \ \ \ \ \ \ \ \ \ \ 
\times
\prod_{i=1}^{(y-k-1)/2}\frac{\left(\frac{x-k}{2}+i\right)_k}{(i)_k}
\prod_{i=1}^{(z-k)/2}\frac{\left(\frac{x+k}{2}+i\right)_{(y-k-1)/2}}{(k+i)_{(y-k-1)/2}}
\frac{\left(\frac{x+k}{2}+i\right)_{(y+k+1)/2}}{(k+i)_{(y+k+1)/2}}.
\label{ebc}
\end{align}
$(${\rm b}$)$. Let $x$, $y$, $z$, $k$ be non-negative integers, $k\leq x,y,z$, and assume that $x$ and $z$ have opposite parity to $y$ and $k$. Then
\begin{align}
&
\M_{\odot}(B'_{x,y,z,k})=P\left(\frac{y+k}{2},\frac{z-k-1}{2},k+1\right)
\nonumber
\\
&\ \ \ \ \ \ \ \ \ \ \ \ 
\times
\prod_{i=0}^{(z-k-1)/2}\frac{\left(\frac{x+k+1}{2}+i\right)_k}{\left(k+i+1\right)_k}
\prod_{i=1}^{(y-k)/2}\frac{\left(\frac{x-k-1}{2}+i\right)_{(z+k+1)/2}}{\left(i\right)_{(z+k+1)/2}}
\frac{\left(\frac{x+k+1}{2}+k+i\right)_{(z-k-1)/2}}{\left(2k+i+1\right)_{(z+k+1)/2}}.
\label{ebd}
\end{align}

\end{theo}

\parindent0pt
\textsc{Remark 4.} Setting $k=0$ in the half of Theorem \ref{tbb} concerning $x$, $y$, $z$, $k$ even, we obtain a product formula for the number of self-complementary plane partitions that fit in a box of even side-lengths $x$, $y$ and $z$. For even $x$ and $z$ and odd $y$, the enumeration of self-complementary plane partitions that fit in this box is obtained from Theorem \ref{tbc}(a) by setting $k=0$, while the case of odd $x$ and $z$ and even $y$ follows by specializing $k=0$ in Theorem \ref{tbc}(b). This amounts to a new proof of the formula enumerating self-complementary plane partitions, first proved by Stanley \cite{Sta}. From this point of view, Theorems \ref{tbb} and \ref{tbc} can be regarded as a generalization of~Stanley's theorem on the number of self-complementary plane partitions that fit in a box.
\parindent15pt

\section{Proof of Theorem 1}

Corresponding to the two sides of equation \eqref{eba}, consider the disjoint union of Cartesian products
\begin{equation}
\mathcal M(G)\times \mathcal M(G\setminus\{a,b,c,d\}) 
\label{eca}
\end{equation}
and
\begin{align}
\mathcal M(G\setminus\{c,d\})\times\mathcal M(G\setminus\{a,b\})
&\cup
\mathcal M(G\setminus\{b,d\})\times\mathcal M(G\setminus\{a,c\})
\nonumber
\\
&\cup
\mathcal M(G\setminus\{a,d\})\times\mathcal M(G\setminus\{b,c\}).
\label{ecb}
\end{align}
For each element $(\mu,\nu)$ of \eqref{eca} or \eqref{ecb}, think of the edges of $\mu$ as being marked with solid lines, and of the edges of $\nu$ as marked by dotted lines, on the same copy of the graph $G$ (the edges that are common to $\mu$ and $\nu$ are marked both solid and dotted, by two parallel arcs). Note that for all $(\mu,\nu)$ corresponding to Cartesian products in \eqref{ecb}, $d$ is matched by a dotted edge (this is the reason for the chosen order of factors in the Cartesian products of \eqref{ecb}).

Define the weight of $(\mu,\nu)$ to be the product of the weight of $\mu$ and the weight of $\nu$. Then the total weight of the elements of the set \eqref{eca} is equal to the left hand side of equation~\eqref{eba}, while the total weight of the elements of the set \eqref{ecb} equals the right hand side of \eqref{eba}. Therefore, to prove \eqref{eba} it suffices to construct a weight-preserving bijection between the sets~\eqref{eca} and~\eqref{ecb}.

We construct such a bijection as follows. Let $(\mu,\nu)$ be an element of \eqref{eca}. Map  $(\mu,\nu)$ to what we get from it by ``shifting along the path containing $d$.'' More precisely, note that when considering the edges of $\mu$ and $\nu$ together on the same copy of $G$, each of the vertices $a,b,c,d$ is incident to precisely one edge. All the other vertices of $G$  are incident to one solid edge and one dotted edge. 

This implies that $\mu\cup\nu$ is the disjoint union of two paths connecting $a$, $b$, $c$ and $d$ among themselves, and cycles covering all the remaining vertices of~$G$; furthermore, both along the paths and along the cycles, the edges alternate between solid and dotted. Consider the path containing $d$, and change each solid edge in it to dotted, and each dotted edge to solid. Denote the resulting pair of matchings by $(\mu',\nu')$. 

Clearly, the weight of $(\mu',\nu')$ is the same as the weight of $(\mu,\nu)$. Therefore, it is enough to show that the map $(\mu,\nu)\mapsto(\mu',\nu')$ is a bijection.

\begin{figure}[h]
  \centerline{
\hfill
{\includegraphics[width=0.80\textwidth]{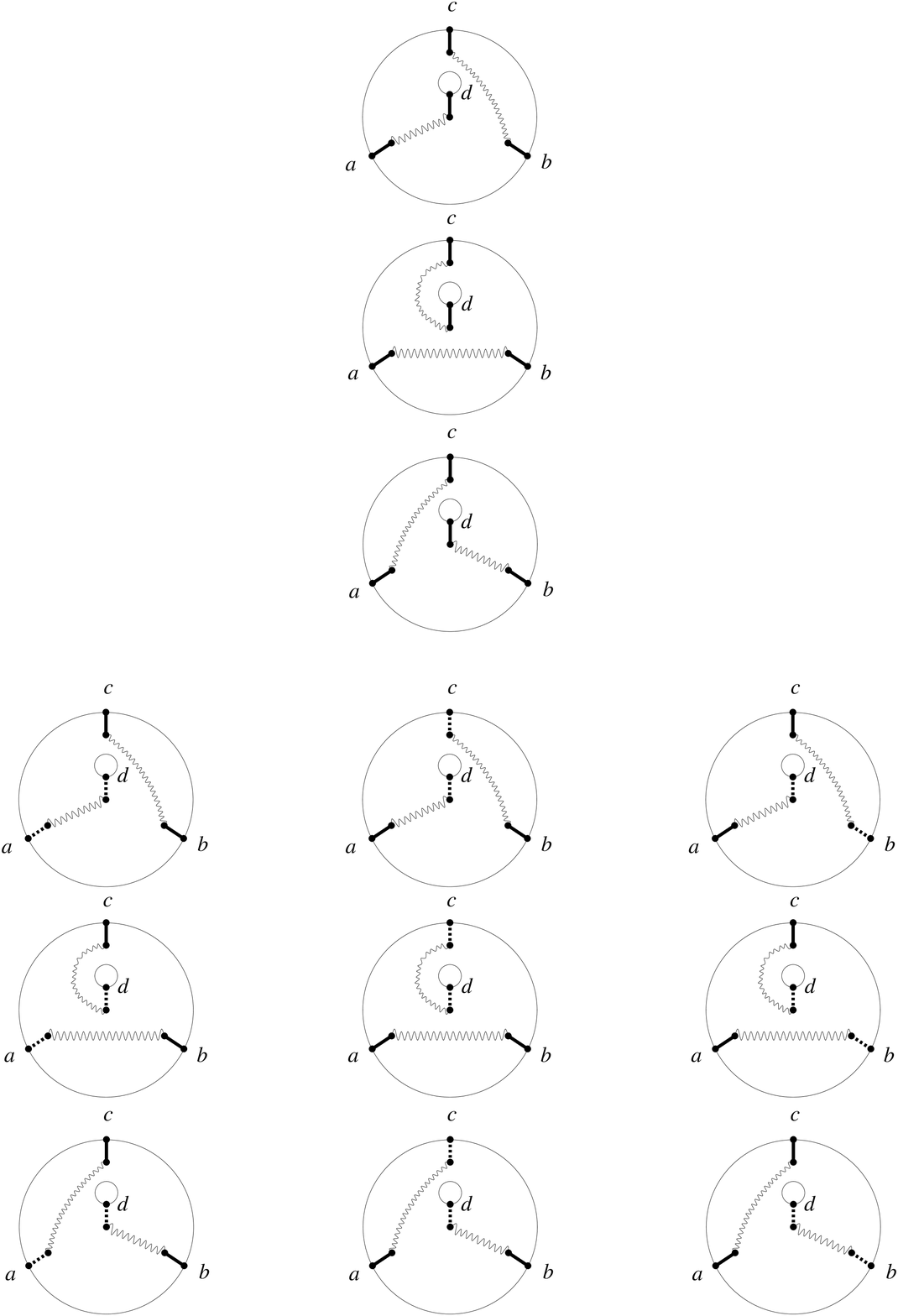}}
\hfill
}
\caption{\label{fca} Schematic representation of the bijection proving Theorem \ref{tba}. Denoting the sets of superpositions of matchings sketched in the top half by $A$, $B$ and $C$ (from top down), and those sketched in the bottom half by  $A'_1$, $B'_1$, $C'_1$ (left column),  $A'_2$, $B'_2$, $C'_2$ (middle column) and  $A'_3$, $B'_3$, $C'_3$ (right column), it turns out that $A$ is in bijection with $A'_1$, $B$ with $B'_2$, and $C$ with $C'_3$, while $A'_2$, $A'_3$, $B'_1$, $B'_3$, $C'_1$ and $C'_2$ are all empty.}
\label{fbb}
\end{figure}

To see this, we partition the Cartesian product \eqref{eca} into three classes, according to the three connection possibilities (in the corresponding superposition of matchings) for vertices $a$, $b$, $c$,~$d$: We gather those $(\mu,\nu)$ for which, in the superposition of $\mu$ and $\nu$, $a$ is connected by a path to $d$,\linebreak into one class; those for which $a$ is connected by a path to $b$ into another class; and those for which $a$ is connected by a path to $c$ into a third class (these partitions are represented schematically in the top half of Figure \ref{fca}).

Partition similarly each of the three Cartesian products in \eqref{ecb} into three classes, according to the same three possible connection types among the vertices $a$, $b$, $c$ and $d$. These are illustrated in the bottom half of Figure \ref{fca}. 

Under the above mapping $(\mu,\nu)\mapsto(\mu',\nu')$, each of the three classes of superpositions of matchings indicated in the top half of Figure \ref{fca} is mapped bijectively to one of the nine classes of superpositions of matchings indicated in the bottom half of Figure \ref{fca}.  The correspondence is indicated in Figure \ref{fca} (the top group of 3 ``balls'' is denoted by $A$, and the bottom groups by $A'$, $B'$ and $C'$; the subscript $i$ indicates that the $i$th ball from the group --- counting from the top --- is chosen).

Indeed, consider for instance the class consisting of those $(\mu,\nu)$ in the Cartesian product in~\eqref{eca} for which $d$ is connected by a path to $a$ in the superposition of $\mu$ and $\nu$ (this is represented by the top ball in the group of three in the top half of Figure \ref{fca}). Since both $d$
%
%
and $a$ are matched by solid edges, after we apply our construction to obtain $(\mu',\nu')$ (which recall consists in reversing the type of all edges in the path containing $d$ --- turning solid edges to dotted, and dotted to solid), both $d$ and $a$ will be matched by dotted edges. They will also clearly still be connected to one another. Therefore this class is mapped into the ``$a$ connected to $d$'' class of the first Cartesian product in \eqref{ecb} (represented by the top ball in the first group of three in the bottom half of Figure \ref{fca}). Since our map is clearly an involution, it establishes a bijection between these two classes.

\parindent15pt
Similarly, the ``$a$ connected to $b$'' class of the Cartesian product in \eqref{eca} is mapped bijectively onto the ``$a$ connected to $b$'' class of the second Cartesian product in \eqref{ecb}, and the ``$a$ connected to $c$'' class of the Cartesian product in \eqref{eca} is mapped bijectively onto the ``$a$ connected to $c$'' class of the third Cartesian product in \eqref{ecb}.

The key to our proof is that, under the hypotheses of our theorem, the remaining six classes of \eqref{ecb} are in fact empty!

To see this, note that in each of these six classes (if they were non-empty!), the path that connects two of the vertices $a$, $b$ and $c$ to one another has an even number of edges (as the starting and ending edges are of opposite type --- one solid, the other dotted ---, and edge types alternate). However, by our assumption, such even length paths cannot arise in superpositions of matchings. Therefore these six classes are indeed empty, and the proof is complete. 

\section{Proof of Theorems \ref{tbaa} and \ref{tbab}}

Note that in both these theorems we may assume that the center of symmetry of $G$ is not at a vertex --- otherwise $G$ would have an odd number of vertices, and all the quantities in the statements would be zero.

For a centrally symmetric planar graph $G$ embedded in a disk, denote by $\tilde{G}$ the orbit graph with respect to the symmetry trough the center. Note that any perfect matching $\tilde{\mu}$ of the orbit graph $\tilde{G}$ can be identified with the centrally symmetric perfect matching $\mu$ of $G$ that maps to $\tilde{\mu}$ under the action of this symmetry.

{\it Proof of Theorem \ref{tbaa}.} Let $\tilde{G}$ be the orbit graph of $G$ under the action of central symmetry, and let $a$, $b$, $c$ and $d$ be the vertices of $\tilde{G}$ corresponding to $\{a_1,a_2\}$, $\{b_1,b_2\}$, $\{c_1,c_2\}$ and $\{d_1,d_2\}$, respectively. The statement of Theorem \ref{tbaa} will follow from Theorem \ref{tba} applied to $\tilde{G}$, provided we show that the condition on the oddness of the path lengths in the hypothesis of Theorem \ref{tba} holds for $\tilde{G}$.

We need to show that, given $S\in\{\{a,b\},\{a,c\},\{b,c\}\}$, for any perfect matchings $\tilde{\mu}\in\mathcal{M}(\tilde{G}\setminus S)$ and $\tilde{\nu}\in\mathcal{M}(\tilde{G}\setminus \overline{S})$, the unique path in the superposition of $\tilde{\mu}$ and $\tilde{\nu}$ that connects two of the vertices $\{a,b,c\}$ to one another has odd length.

For definiteness, suppose $S=\{a,b\}$. Viewing $\tilde{\mu}$ and $\tilde{\nu}$ as centrally symmetric perfect matchings $\mu$ and $\nu$ of $G\setminus\{a_1,a_2,b_1,b_2\}$ and $G\setminus\{c_1,c_2,d_1,d_2\}$, respectively, this amounts to proving that in the superposition of $\mu$ and $\nu$, the unique pair of paths $(P_1,P_2)$ (with $P_2$ the image of $P_1$ through the center) that connect vertices of $\{a_1,b_1,c_1,a_2,b_2,c_2\}$ among themselves have odd length.

\begin{figure}[h]
  \centerline{
\hfill
{\includegraphics[width=0.23\textwidth]{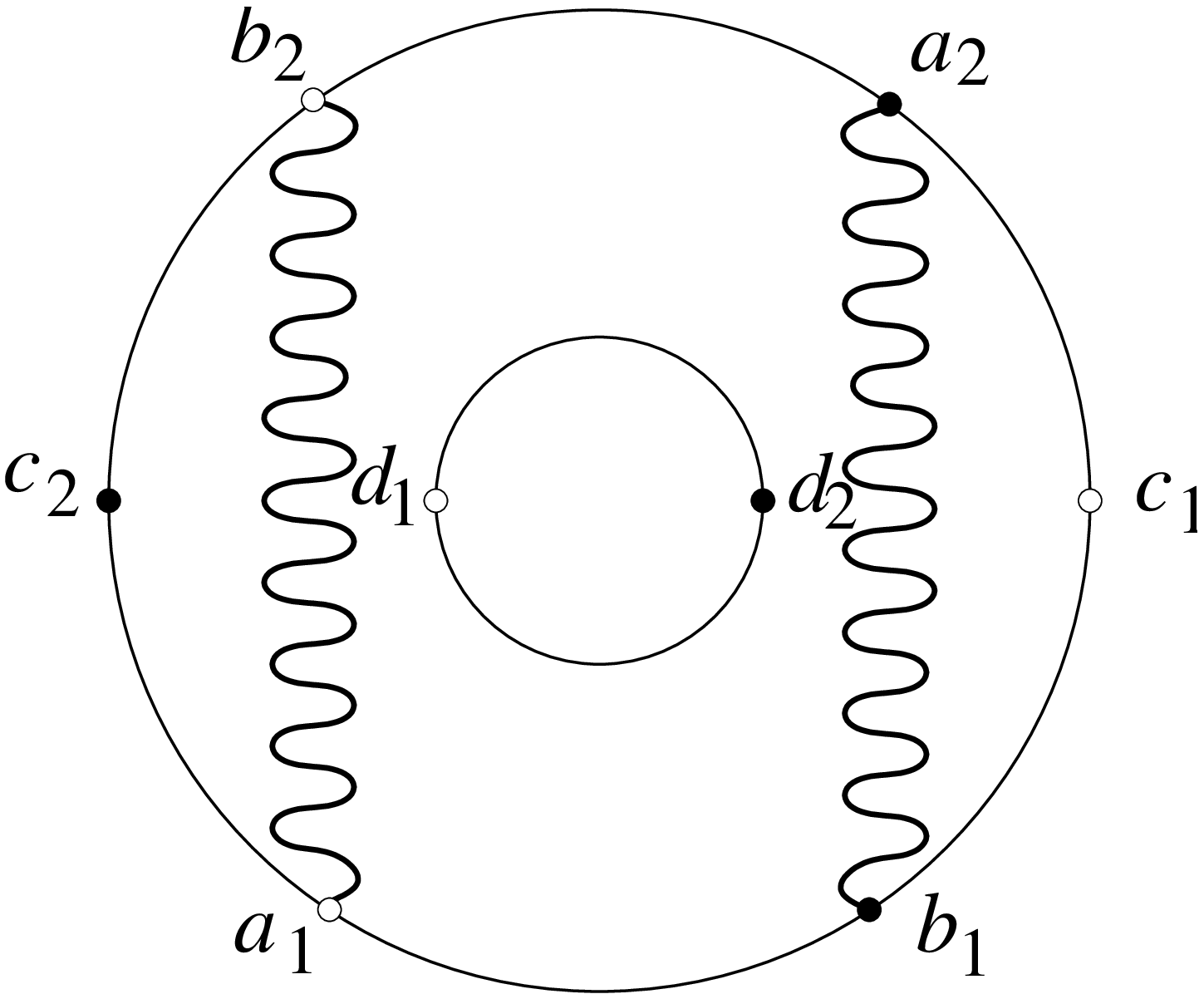}}
\hfill
{\includegraphics[width=0.23\textwidth]{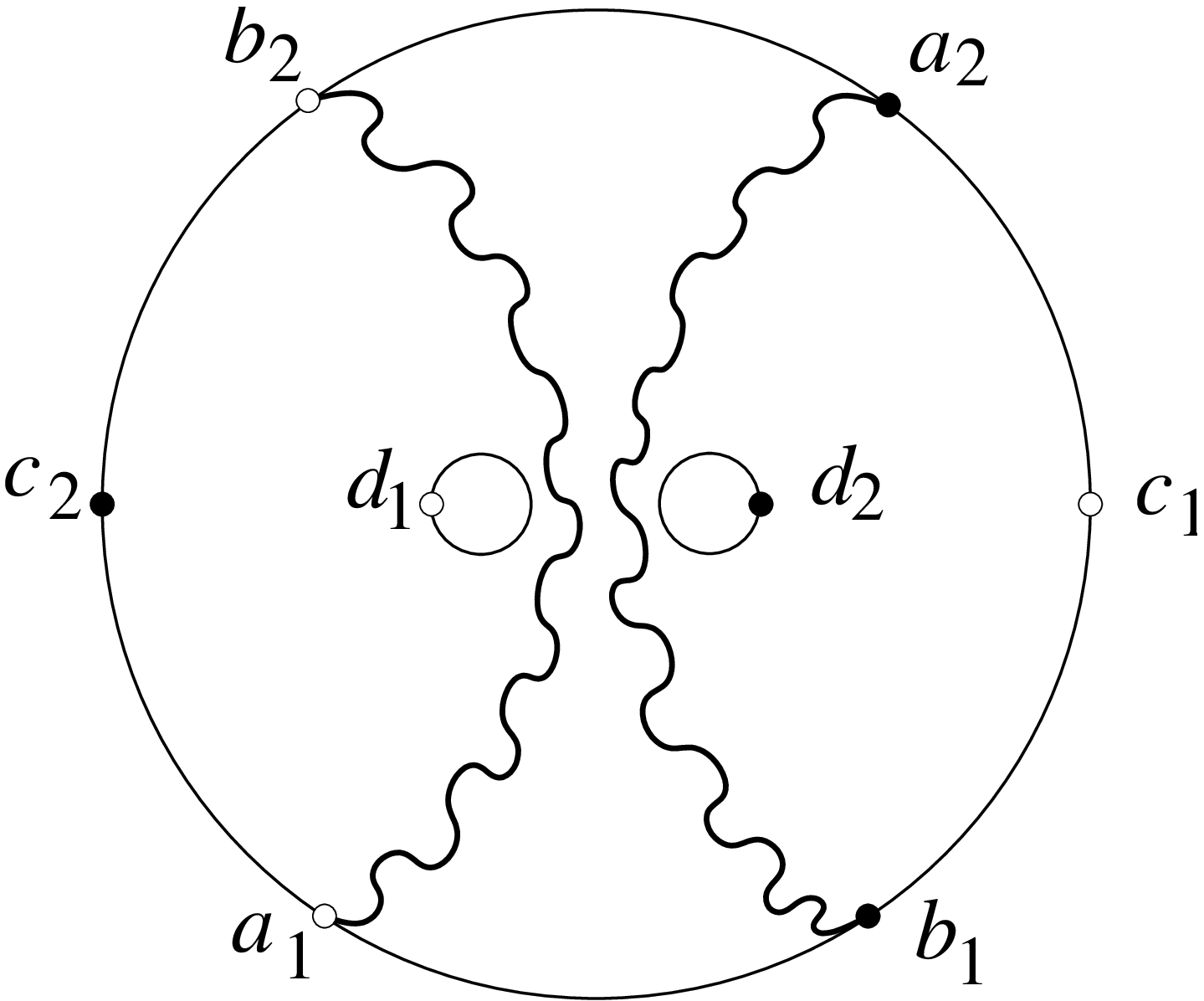}}
\hfill
}
  \caption{\label{fda} A connection type that does not occur from superpositions of matchings (left). A connection type arising from the superposition of a special pair of matchings (right).}
\end{figure}

Suppose towards a contradiction that $P_1$ has even length. Since $G$ is bipartite, the endpoints of $P_1$ must have the same color. Because the colors of the vertices $\{a_1,b_1,c_1,a_2,b_2,c_2\}$ alternate, $P_1$ must connect two of these vertices that are next nearest neighbors in this cyclic order (see the picture on the left in Figure \ref{fda}). Since the image $P_2$ of $P_1$ through the center must be disjoint from $P_1$, it follows that $P_1$ must pass between the inner face and the unique member $v$ of  $\{a_1,b_1,c_1,a_2,b_2,c_2\}$ in between the endpoints of $P_1$. But this would prevent $v$ from being connected to $d_1$ or $d_2$ in $\mu\cup\nu$ by a path disjoint from $P_1$, a contradiction. \epf

{\it Proof of Theorem \ref{tbab}.} Identify the centrally symmetric perfect matchings of $G\setminus S$, $S\subset\{a,b,c,d\}$, $|S|$ even, with the corresponding matchings of the orbit graph $\tilde{G}$.

In the proof of Theorem \ref{tba}, we partitioned each of the four Cartesian products in \eqref{eca} and \eqref{ecb} into three subclasses, depending on the way vertices $a$, $b$, $c$ and $d$ are connected by paths in the superposition of the corresponding matchings. We saw that each of the three subclasses corresponding to \eqref{eca} is in bijection with one subclass of the three Cartesian products in \eqref{ecb}. The proof followed because, under the hypothesis of Theorem \ref{tba}, the remaining six subclasses were empty.


In general, these six classes are not empty, and therefore contribute to the total weight of the set of pairs of matchings \eqref{ecb}. We claim that the quantity they contribute is precisely the term that is subtracted on the right hand side of \eqref{ebab}. Once this claim is proved, the statement of Theorem~\ref{tbab} will follow from the proof of Theorem \ref{tba}.

To prove the claim, note that what characterizes the six classes mentioned above is that in each of them, the length of the unique path in the corresponding superposition of matchings $\tilde{\mu}\cup\tilde{\nu}$ which connects two vertices of $\{a,b,c\}$ with one another is even. This implies that in the superposition of the corresponding centrally symmetric matchings $\mu$ and $\nu$ of $G$, the unique pair of paths $(P_1,P_2)$ (with $P_2$ the image of $P_1$ through the center) that connect vertices of $\{a_1,b_1,c_1,a_2,b_2,c_2\}$ among themselves have even length, and therefore connect vertices of the same color in $G$ (as $G$ is bipartite).

In the superposition $\mu\cup\nu$ of $\mu$ and $\nu$, the eight vertices $a_1,a_2,b_1,b_2,c_1,c_2,d_1,d_2$ are connected up by four disjoint paths, two of which are $P_1$ and $P_2$, connecting up four of the vertices $\{a_1,b_1,c_1,a_2,b_2,c_2\}$. Since $P_1$ and $P_2$ have same-color endpoints, one readily sees that in order for the remaining four vertices to be connectable by two more paths $P_3$ and $P_4$ so that $P_1$, $P_2$, $P_3$ and $P_4$ are disjoint, $P_1$ and $P_2$ must pass through the faces $F_1$ and $F_2$ which contain $d_1$ and $d_2$ (see the picture on the right in Figure \ref{fda}). This proved the claim, and the proof of Theorem~\ref{tbab} is complete. \epf

\section{Proof of Theorems \ref{tbb} and \ref{tbc}}

Our proofs of Theorems \ref{tbb} and \ref{tbc} are by induction on the sum $x+y+z$ of the side-lengths of the outer hexagon in the regions $B_{x,y,z,k}$ and $B'_{x,y,z,k}$. We identify lattice regions with their planar duals, and lozenge tilings of the former with perfect matchings of the latter.

For the regions $B_{x,y,z,k}$, at the induction step we use the recurrence provided by Theorem~\ref{tbaa}.
Despite the last term in \eqref{ebab}, Theorem \ref{tbab} turns out to give a recurrence for the regions $B'_{x,y,z,k}$. This is because the removed lobes are only distance 1 apart, and therefore there is no room for two disjoint paths to pass between them. Thus the last term on the right hand side of \eqref{ebab} is zero; the resulting recurrence for the number of perfect matchings of the regions $B'_{x,y,z,k}$ is what we use at the induction step in our proof of Theorem \ref{tbc}.

\begin{figure}[h]
  \centerline{
\hfill
{\includegraphics[width=0.37\textwidth]{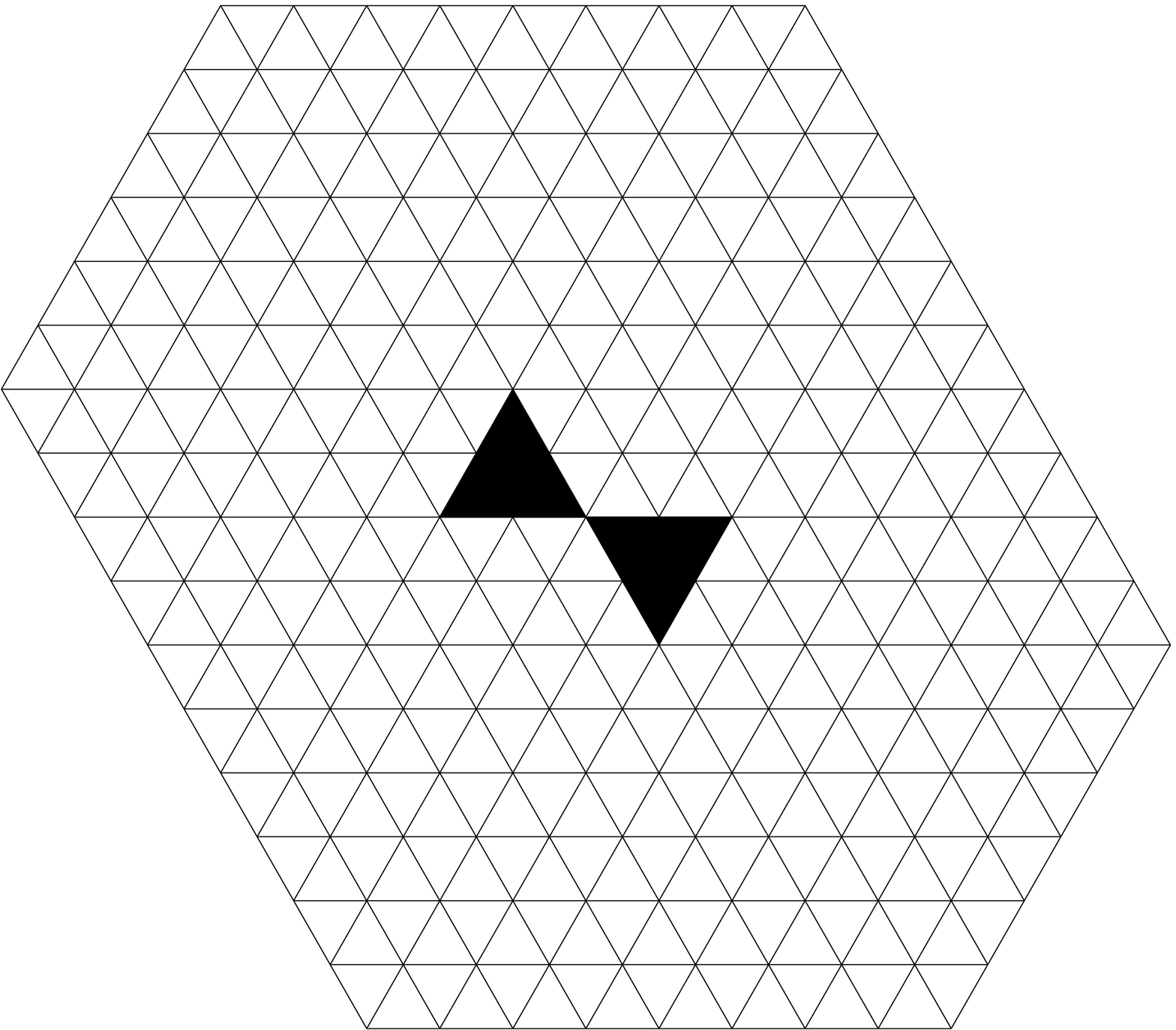}}
\hfill
{\includegraphics[width=0.37\textwidth]{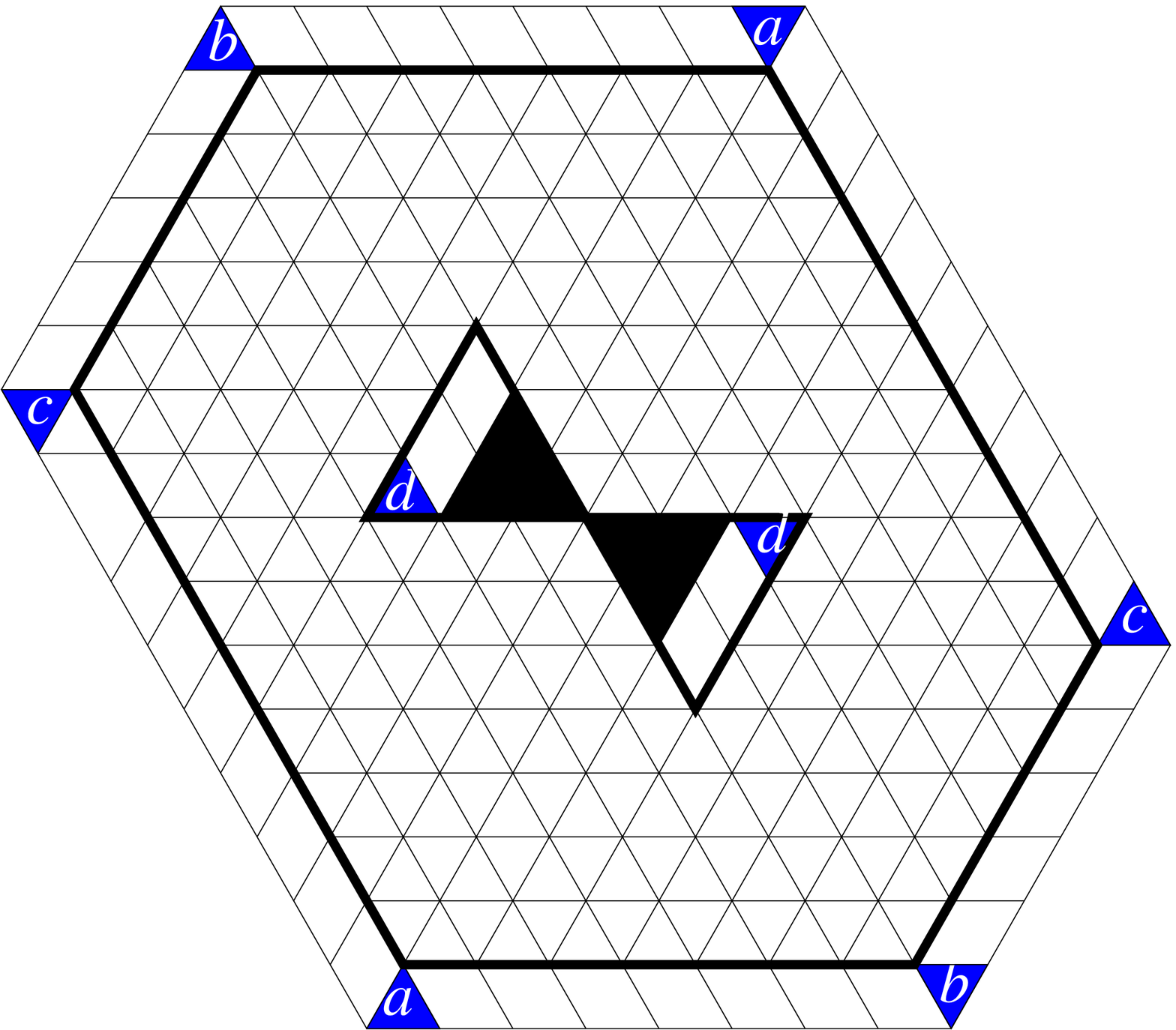}}
\hfill
}
\vskip0.06in
  \centerline{
\hfill
{\includegraphics[width=0.37\textwidth]{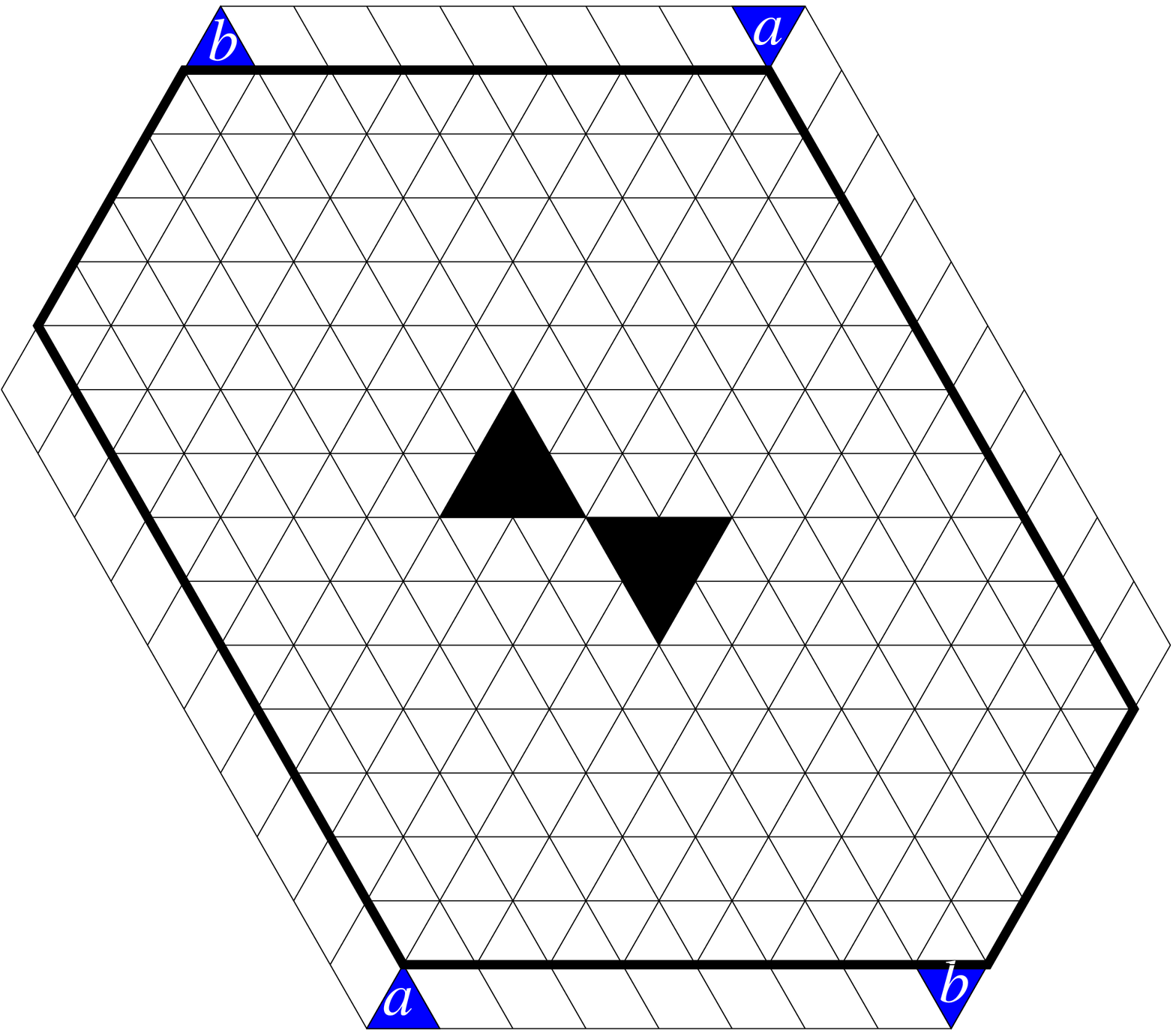}}
\hfill
{\includegraphics[width=0.37\textwidth]{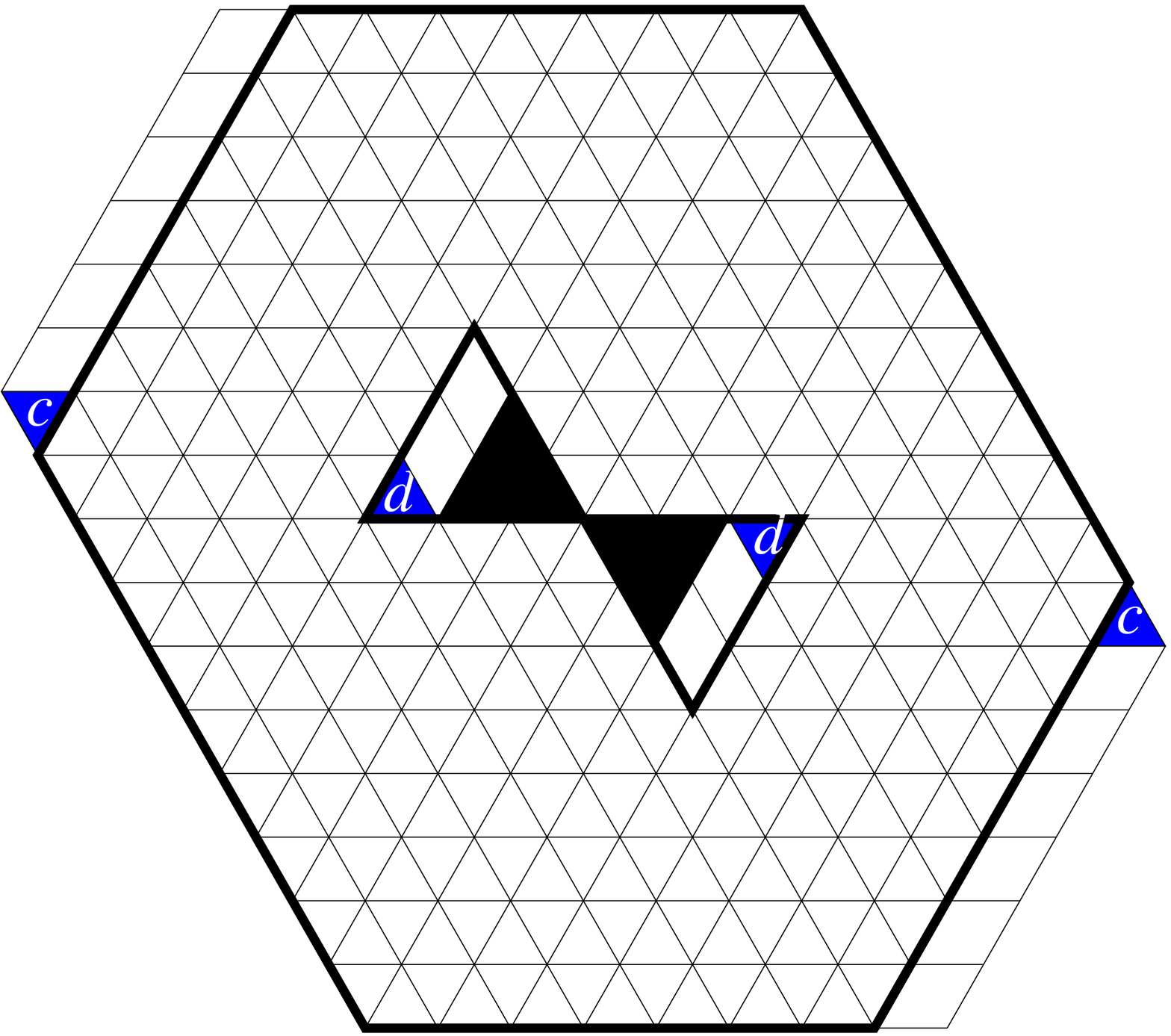}}
\hfill
}
\vskip0.06in
  \centerline{
\hfill
{\includegraphics[width=0.37\textwidth]{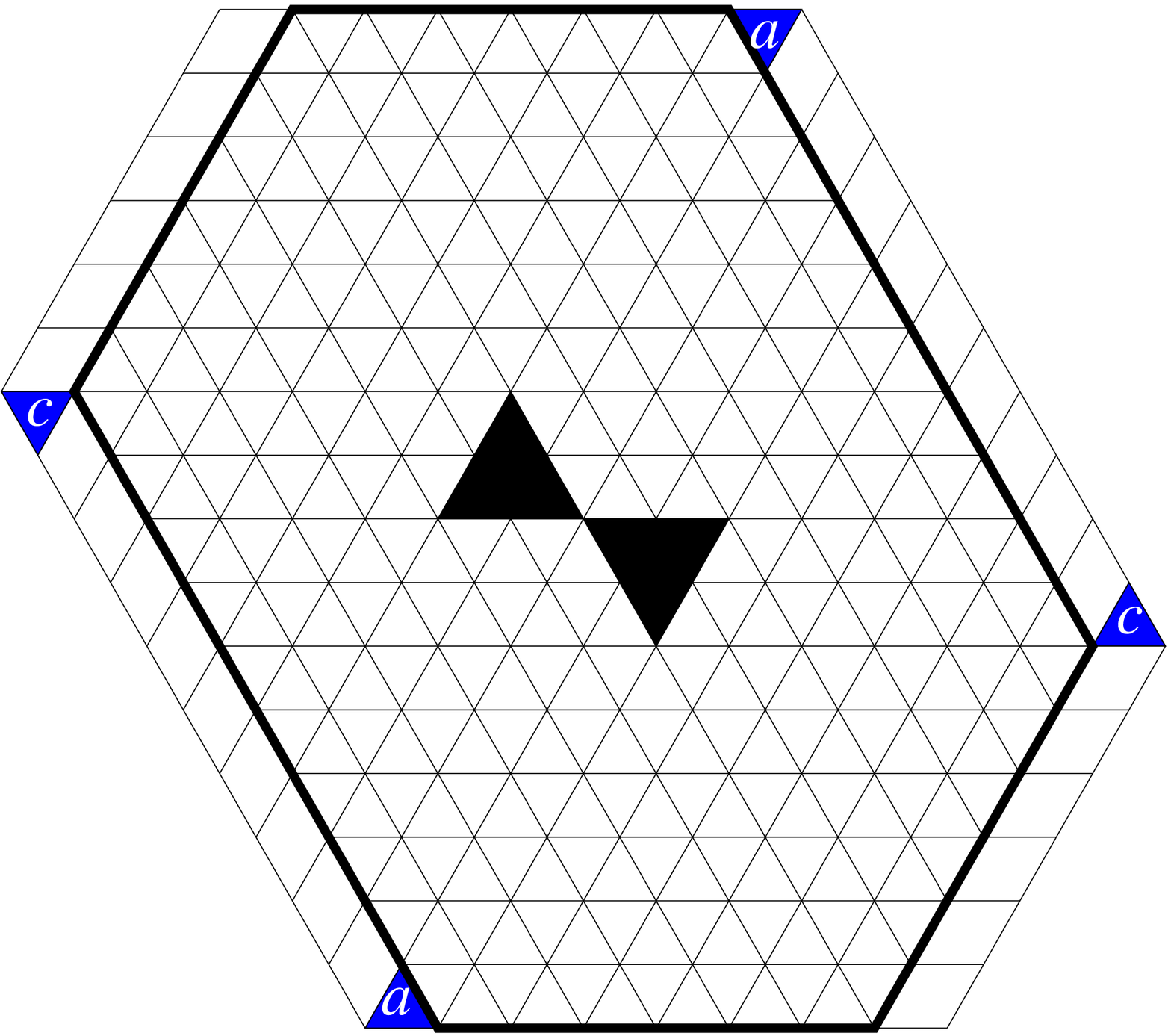}}
\hfill
{\includegraphics[width=0.37\textwidth]{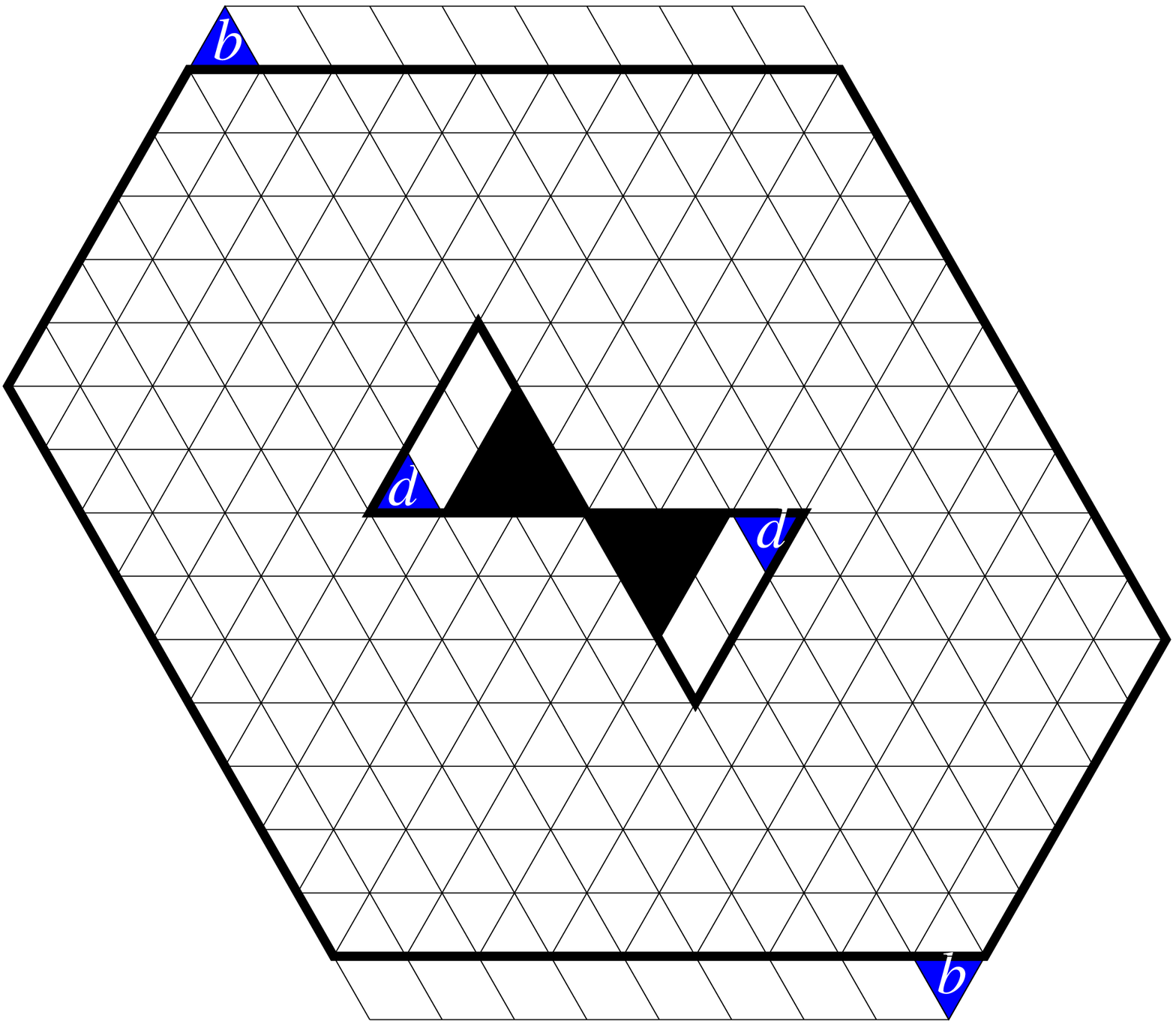}}
\hfill
}
\vskip0.06in
  \centerline{
\hfill
{\includegraphics[width=0.37\textwidth]{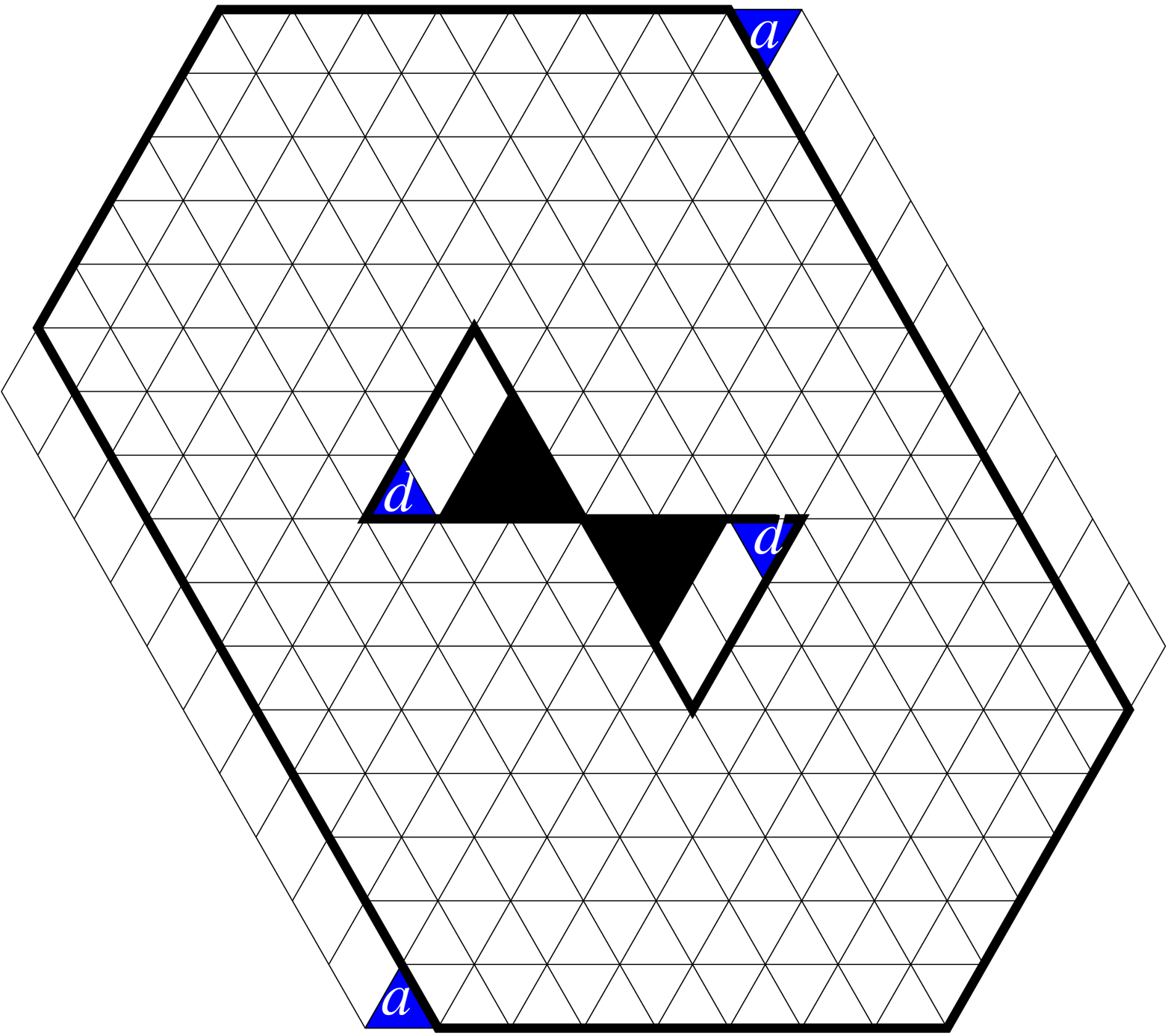}}
\hfill
{\includegraphics[width=0.37\textwidth]{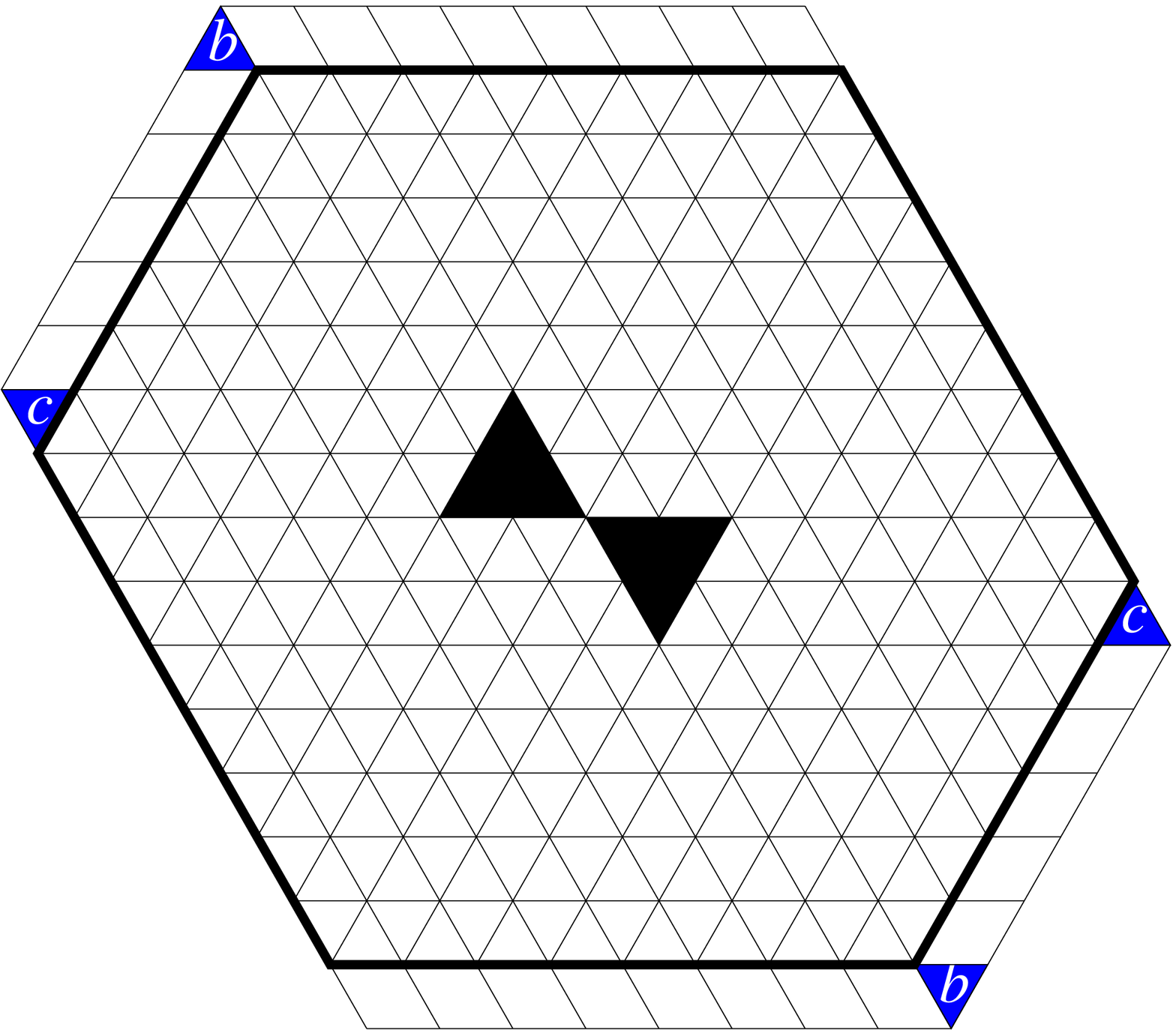}}
\hfill
}
  \caption{\label{fea} Obtaining the recurrence for the $B_{x,y,z,k}$ regions.}
\end{figure}

{\it Proof of Theorem \ref{tbb}}. Consider the region $B_{x,y,z,k}$, and choose the unit triangles $a_1$, $b_1$, $c_1$, $a_2$, $b_2$, $c_2$, $d_1$ and $d_2$ as shown in Figure \ref{fea} (in the figure the indices are dropped; $a_1$, $b_1$ and $c_1$ are the bottom left, bottom right, and rightmost of the marked unit triangles, respectively). Then the hypotheses of Theorem \ref{tbaa} are met, so \eqref{ebaa} holds. After removing the forced lozenges, it follows from Figure \ref{fea} that \eqref{ebaa} takes the form
%
\begin{align}
\M_{\odot}(B_{x,y,z,k})\M_{\odot}(B_{x-1,y-1,z-1,k+1})=
&\M(B_{x,y,z-2,k})\M(B_{x-1,y-1,z+1,k+1})+
\nonumber
\\[5pt]
&\M_{\odot}(B_{x-2,y,z,k})\M_{\odot}(B_{x+1,y-1,z-1,k+1})+
\nonumber
\\[5pt]
&
\M_{\odot}(B_{x-1,y+1,z-1,k+1})\M_{\odot}(B_{x,y-2,z,k}).
\label{eea}
\end{align}

\begin{figure}[h]
  \centerline{
\hfill
{\includegraphics[width=0.30\textwidth]{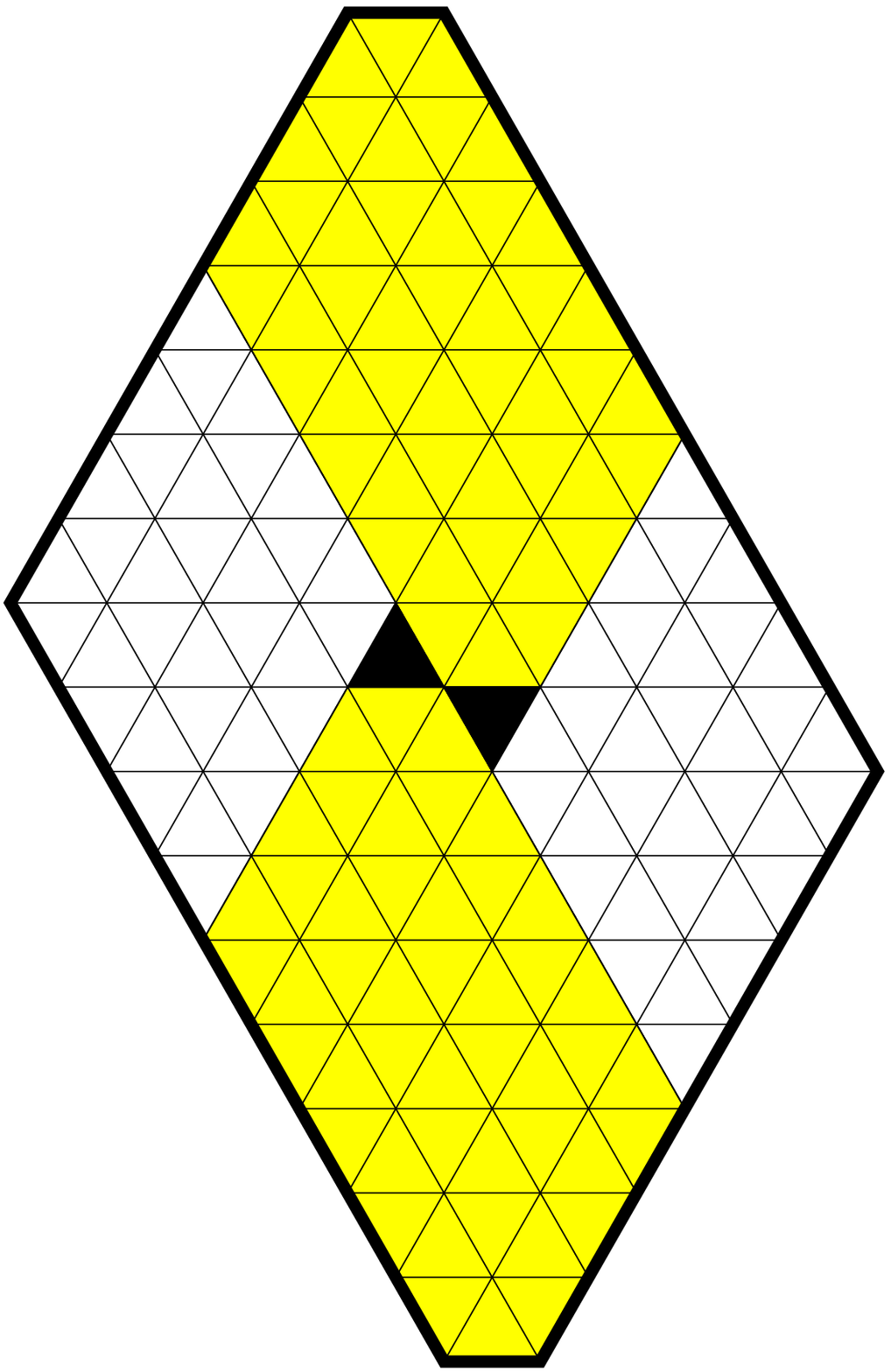}}
\hfill
{\includegraphics[width=0.44\textwidth]{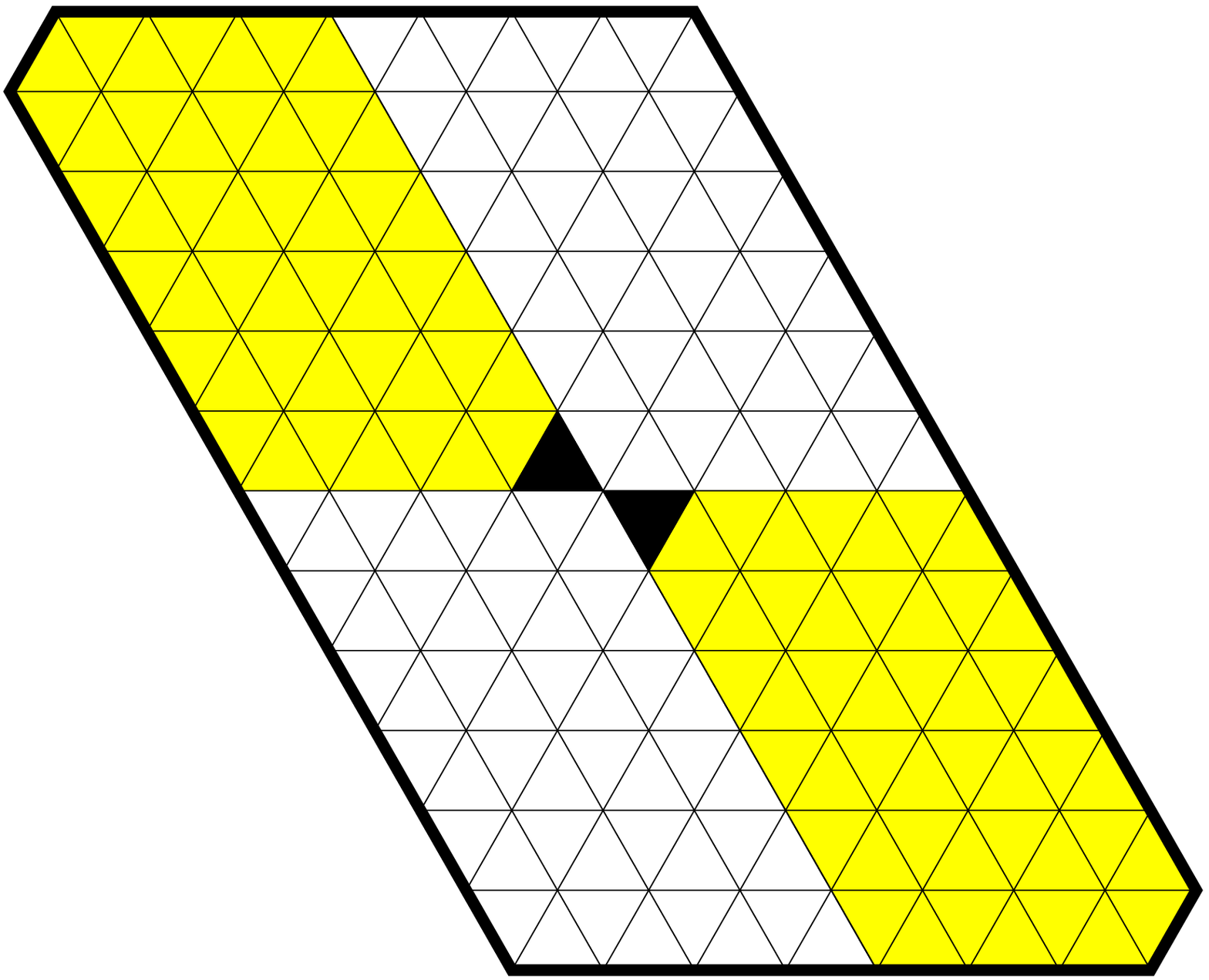}}
\hfill
}
  \caption{\label{feb} The base cases $x=1$ (left) and $z=1$ (right).}
\end{figure}

Note that in \eqref{eea}, the sum of the $x$-, $y$- and $z$-parameters of the first region is $x+y+z$, while for each of the remaining seven regions it is strictly less than $x+y+z$. We will regard \eqref{eea} as providing a recurrence relation for $\M(B_{x,y,z,k})$ in terms of smaller (as measured by the sum of the $x$-, $y$- and $z$-parameters) $B$-regions, and prove \eqref{ebaa} by induction on $x+y+z$.

In order for all the regions in \eqref{eea} to be defined, we need to have $x,y,z\geq2$, and also $k\leq x-2,y-2,z-2$. 
Therefore, the cases when $x$, $y$ or $z$ are equal to 0 or 1, and also when $k$ is equal to one of $x$, $y$ and $z$ ($k$ cannot be one unit less that $x$, $y$ or $z$, as $x,y,z,k$ have the same parity), will be base cases of our induction. Note also that for the seven regions involved in \eqref{eea} besides $B_{x,y,z,k}$, the sum of their $x$-, $y$- and $z$-parameters can be one, two or three units less than $x+y+z$. Thus the cases when $x+y+z$ equals 0, 1 or 2 need to be base cases as well. But clearly, if $0\leq x+y+z\leq2$, then at least one of $x$, $y$ and $z$ must be 0, so in fact the former base cases also cover the latter.

Due to the direction of the removed bowtie, parameters $x$ and $y$ play symmetrical roles. Therefore, it is enough to consider the base cases $x=0$, $x=1$, $z=0$ and $z=1$.

Since $k\leq x,y,z$, if $x=0$, then we must have $k=0$, and the region becomes $B_{0,y,z,0}$, which is a parallelogram, and thus has a unique lozenge tiling. The expression on the right hand side of \eqref{ebb} is also equal to 1 for $x=k=0$, and therefore \eqref{ebb} holds in this case. The case $z=0$ is verified the same way.

Suppose now that $x=1$. Since $k$ must have the same parity as $x$,
it follows that $k=1$. It is easy to see that the shaded hexagons on the left in Figure \ref{feb} are internally tiled in any lozenge tiling of $B_{1,y,z,1}$. Since on the rest of the region the tiling is forced, we obtain
\begin{equation}
\M_{\odot}(B_{1,y,z,1})=\M\left(H_{1,\frac{y+1}{2},\frac{z-1}{2}}\right),
\end{equation}  
and equality  \eqref{ebb} follows from MacMahon's theorem \eqref{eaa}.


For $z=1$, by the same argument as in the previous paragraph we must have $k=1$. The shaded hexagons indicated on the right in Figure \ref{feb} must be internally tiled. Again, the tiling is forced on the remaining part of the region, yielding
\begin{equation}
\M_{\odot}(B_{1,y,z,1})=\M\left(H_{1,\frac{y-1}{2},\frac{z-1}{2}}\right),
\end{equation}  
and equality  \eqref{ebb} follows again from MacMahon's theorem \eqref{eaa}.

The cases when $k$ is equal to $x$, $y$ or $z$ follow in a similar fashion, as they imply that two symmetric hexagons with one pair of opposite sides of length $k$ are internally tiled.

For the induction step, assume that equality \eqref{ebb} holds for all non-negative integers $x,y,z\geq2$, $k\leq x,y,z$ of the same parity, for which the sum of the $x$-, $y$- and $z$-parameters is strictly less than $x+y+z$. We need to prove that it also holds for $B_{x,y,z,k}$.

Use equation \eqref{eea} to express $\M(B_{x,y,z,k})$ in terms of the other seven $B$-regions involved. The induction hypothesis can be applied to each of them, yielding a concrete simple product for the number of its lozenge tilings. It is routine to verify that the resulting expression for $\M(B_{x,y,z,k})$ agrees with the right hand side of \eqref{ebb}. This completes the induction step. \epf

\begin{figure}[h]
  \centerline{
\hfill
{\includegraphics[width=0.32\textwidth]{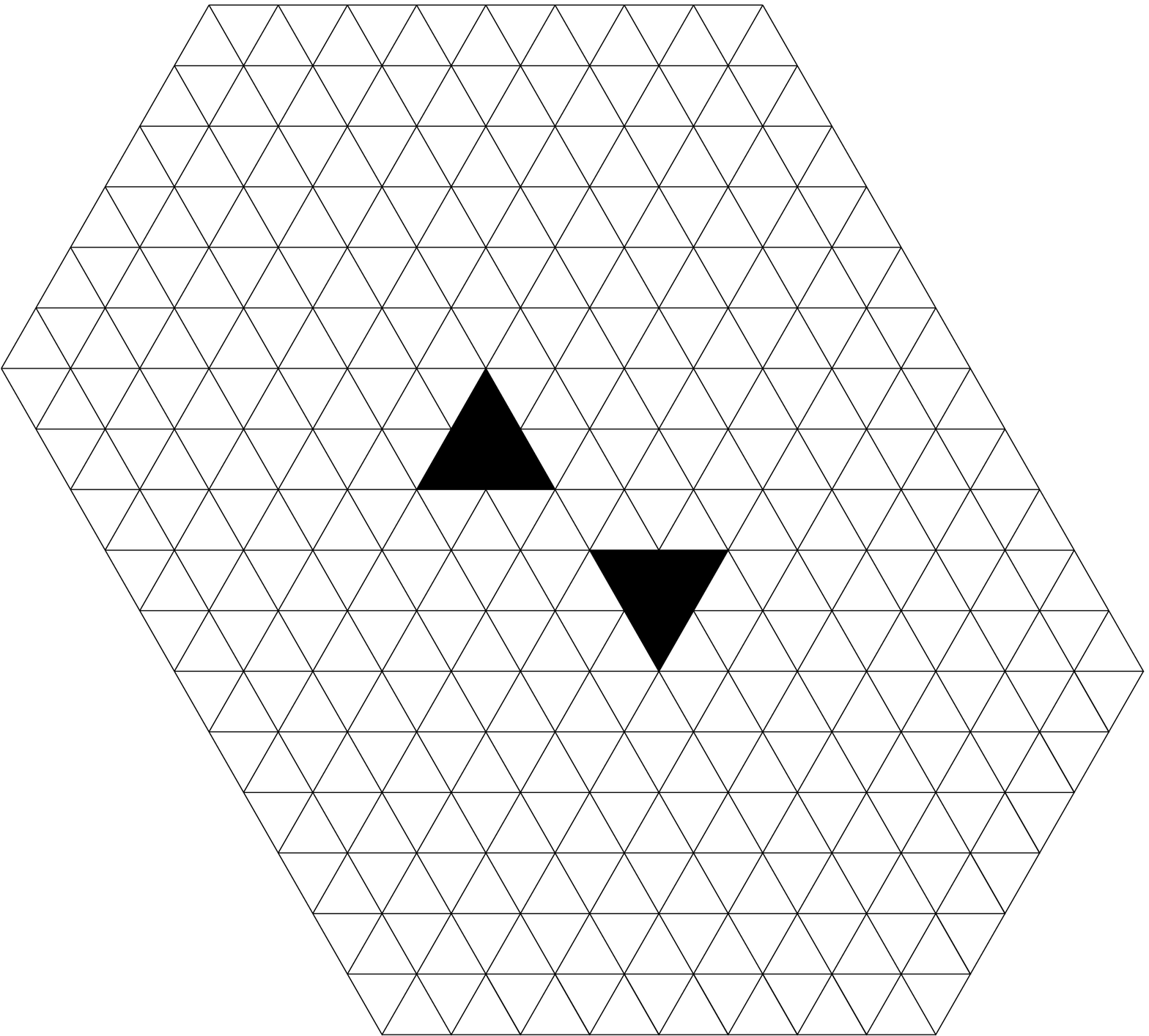}}
\hfill
{\includegraphics[width=0.32\textwidth]{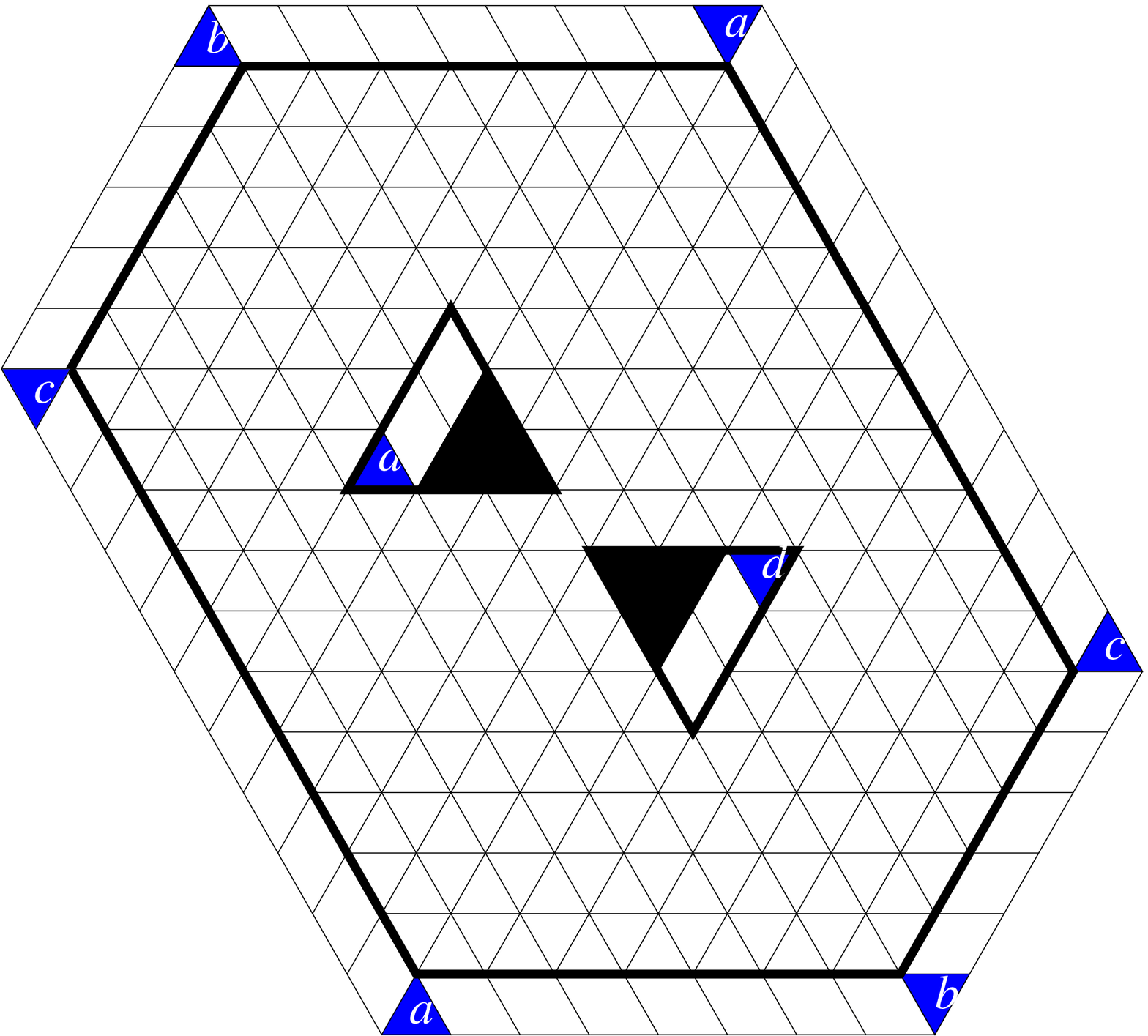}}
\hfill
}
\vskip0.05in
  \centerline{
\hfill
{\includegraphics[width=0.32\textwidth]{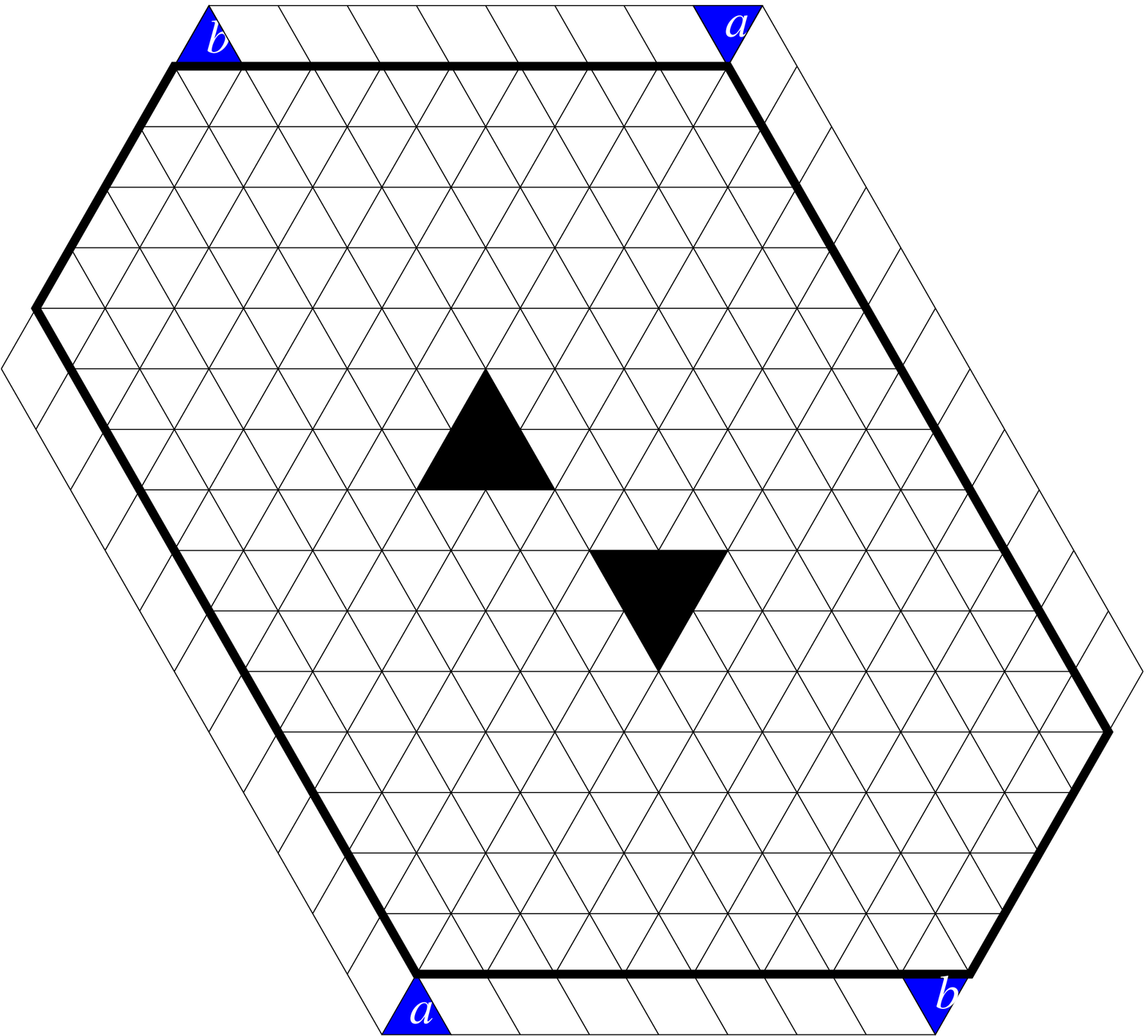}}
\hfill
{\includegraphics[width=0.32\textwidth]{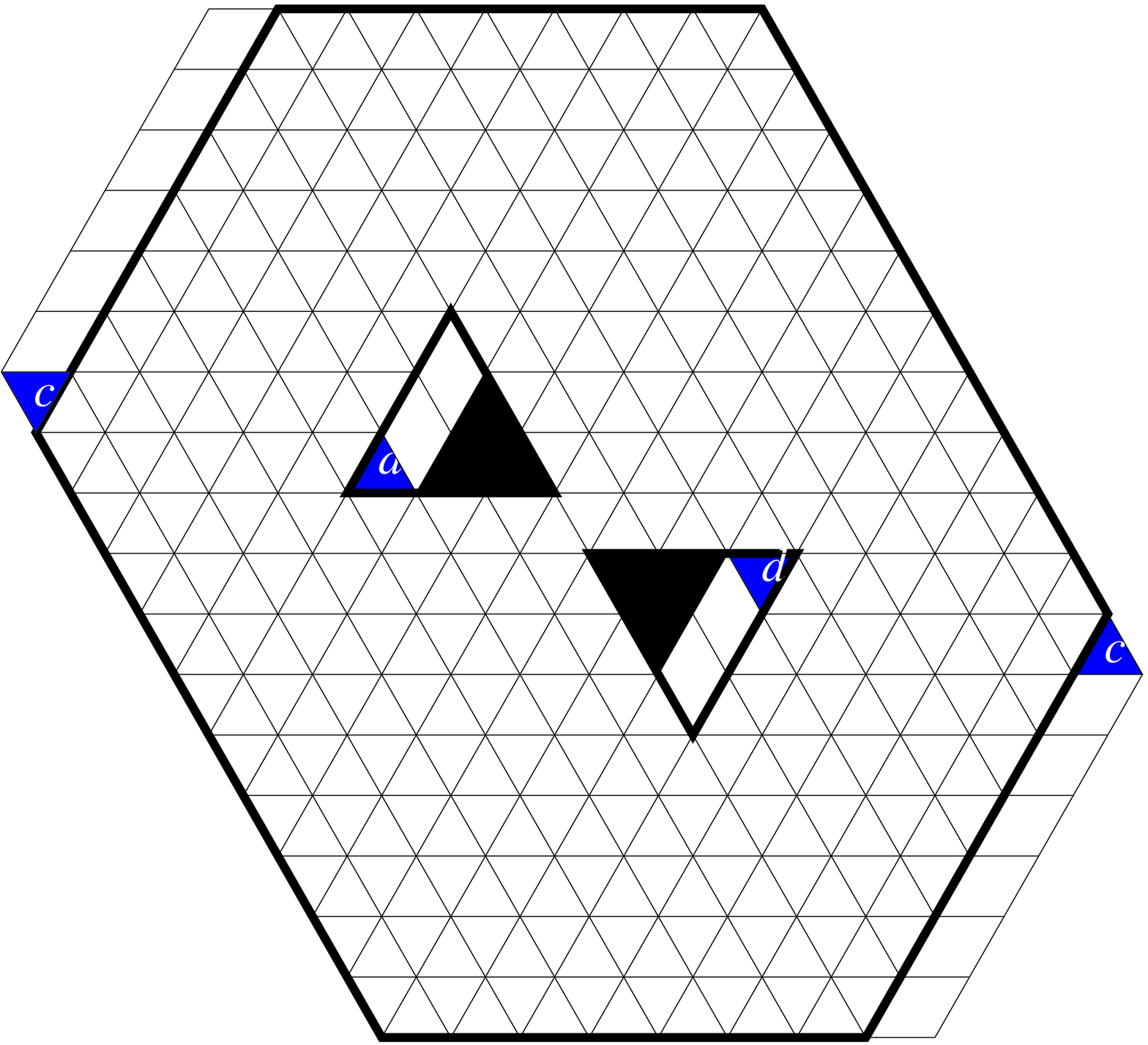}}
\hfill
}
\vskip0.05in
  \centerline{
\hfill
{\includegraphics[width=0.32\textwidth]{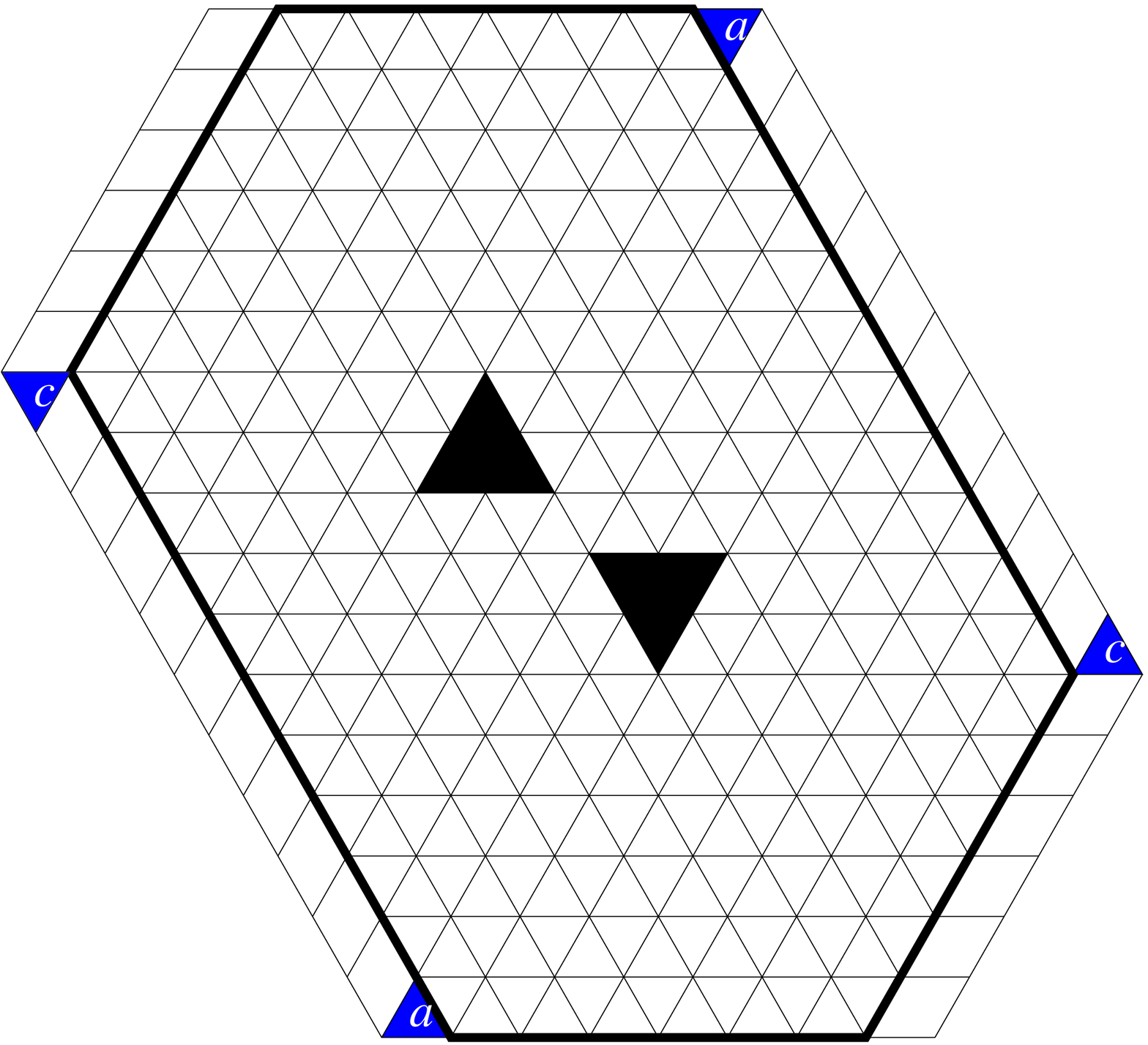}}
\hfill
{\includegraphics[width=0.32\textwidth]{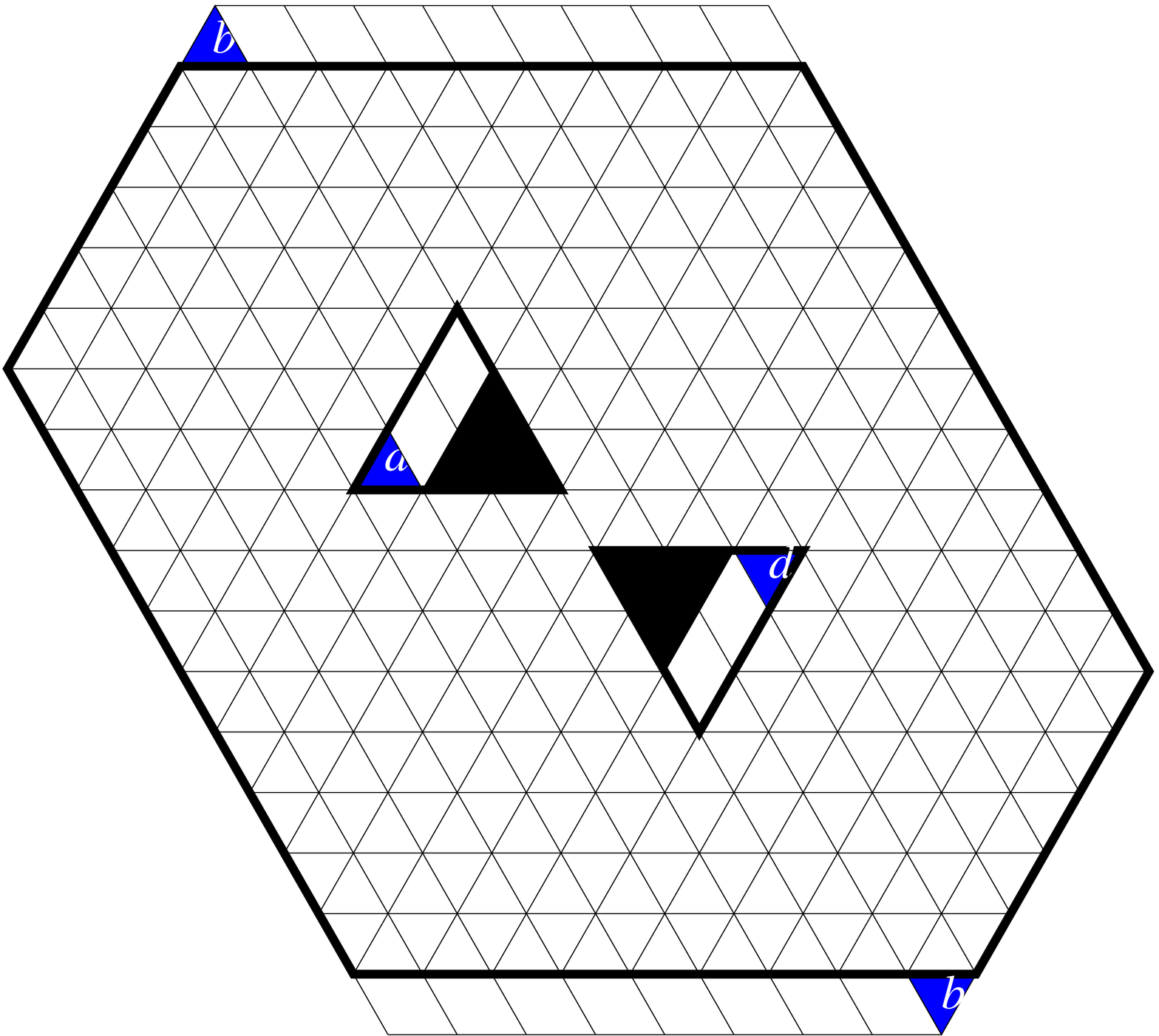}}
\hfill
}
\vskip0.05in
  \centerline{
\hfill
{\includegraphics[width=0.32\textwidth]{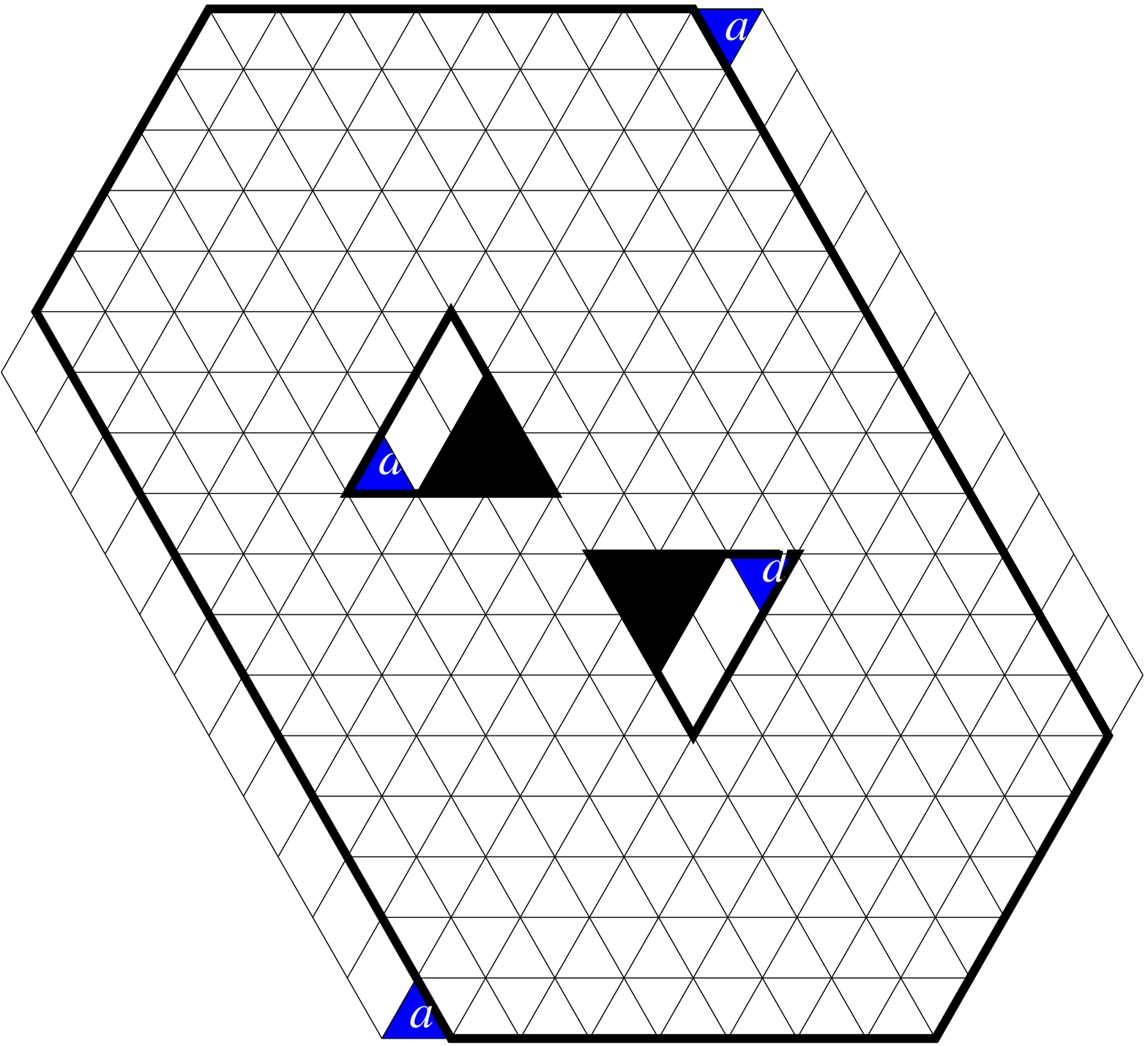}}
\hfill
{\includegraphics[width=0.32\textwidth]{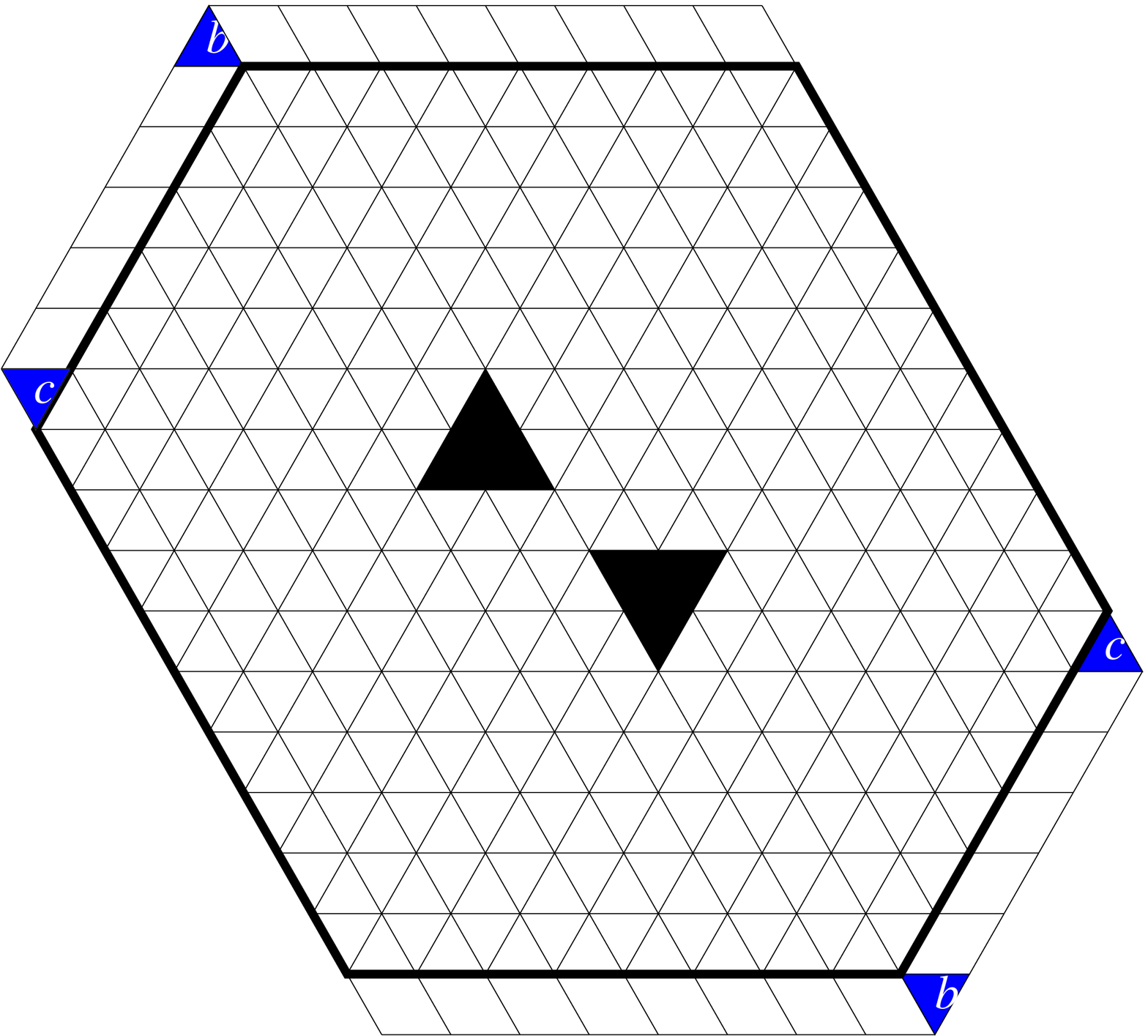}}
\hfill
}
  \caption{\label{fec} Obtaining the recurrence for the $B'_{x,y,z,k}$ regions.}
\end{figure}

{\it Proof of Theorem \ref{tbc}}. Let $G$ be the (planar) dual graph of the region $B'_{x,y,z,k}$, and choose the vertices $a_1$, $b_1$, $c_1$, $a_2$, $b_2$, $c_2$, $d_1$ and $d_2$ of $G$ to correspond to the unit triangles of $B'_{x,y,z,k}$ indicated in Figure \ref{fec} (as before, in the figure the indices are dropped; $a_1$, $b_1$ and $c_1$ are the bottom left, bottom right, and rightmost of the marked unit triangles, respectively). Then the hypotheses of Theorem \ref{tbab} are met. Furthermore, since there is no room in between the two halves of the removed disconnected bowtie for two disjoint paths to pass through, it follows that the quantity that is subtracted on the right hand side of \eqref{tbab} is equal to zero. Therefore, Theorem \ref{tbab} gives
\\[1pt]
\begin{align}
\M_{\odot}(B'_{x,y,z,k})\M_{\odot}(B'_{x-1,y-1,z-1,k+1})=
&\M(B'_{x,y,z-2,k})\M(B'_{x-1,y-1,z+1,k+1})+
\nonumber
\\[5pt]
&\M_{\odot}(B'_{x-2,y,z,k})\M_{\odot}(B'_{x+1,y-1,z-1,k+1})+
\nonumber
\\[5pt]
&\M_{\odot}(B'_{x-1,y+1,z-1,k+1})\M_{\odot}(B'_{x,y-2,z,k}),
\label{eed}
\end{align}
which is precisely the same as the recurrence \eqref{eea} satisfied by the number of lozenge tilings of the regions $B_{x,y,z,k}$. By the considerations at the beginning of the proof of Theorem \ref{tbaa}, formulas \eqref{ebc} and  \eqref{ebd} can be proved by induction on $x+y+z$. The base cases are the cases when $x$, $y$ or $z$ are equal to 0 or 1, together with the cases $k=x$, $k=z$ and $k=y-1$ for part a), and $k=x-1$, $k=z-1$ and $k=y$ for part b).


We discuss below in detail the base cases $x=0$, $x=1$, $y=0$ and $y=1$. The cases $z=0$ and $z=1$ follow using the same arguments.

The base cases $x=0$ and $y=0$ follow the same way as the corresponding base cases in the proof of Theorem \ref{tbaa}, because since $k\leq x,y$, we must have $k=0$ in these cases.

However, there are some changes for the base cases $x=1$ and $y=1$, because now the lobes of the removed bowtie are not touching, and also because $k$ is now free to have any parity.

Suppose $x=1$. Then since $k\leq x$, $k$ is either 0 or 1. If $k=0$, the removed disconnected bowtie becomes empty. Also, since in this case $x$ and $k$ have opposite parities, we are in the situation covered by part (b) of the theorem. The region $B'_{1,y,z,0}$ is the hexagon $H_{1,y,z}$, and we need to verify that its number of centrally symmetric lozenge tilings is given by the $x=1$, $k=0$ specialization of formula \eqref{ebd}.

\begin{figure}[h]
  \centerline{
\hfill
{\includegraphics[width=0.30\textwidth]{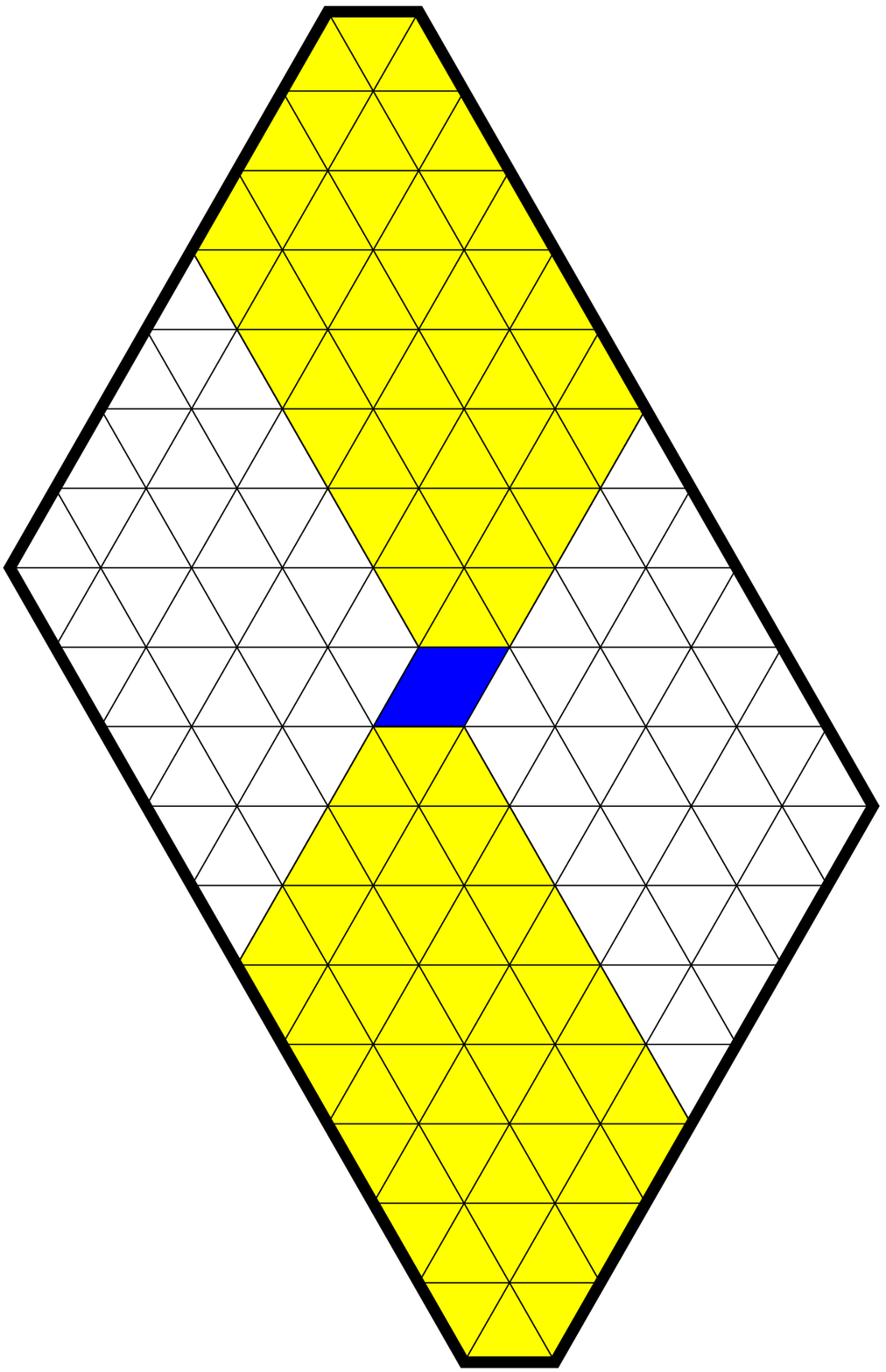}}
\hfill
{\includegraphics[width=0.30\textwidth]{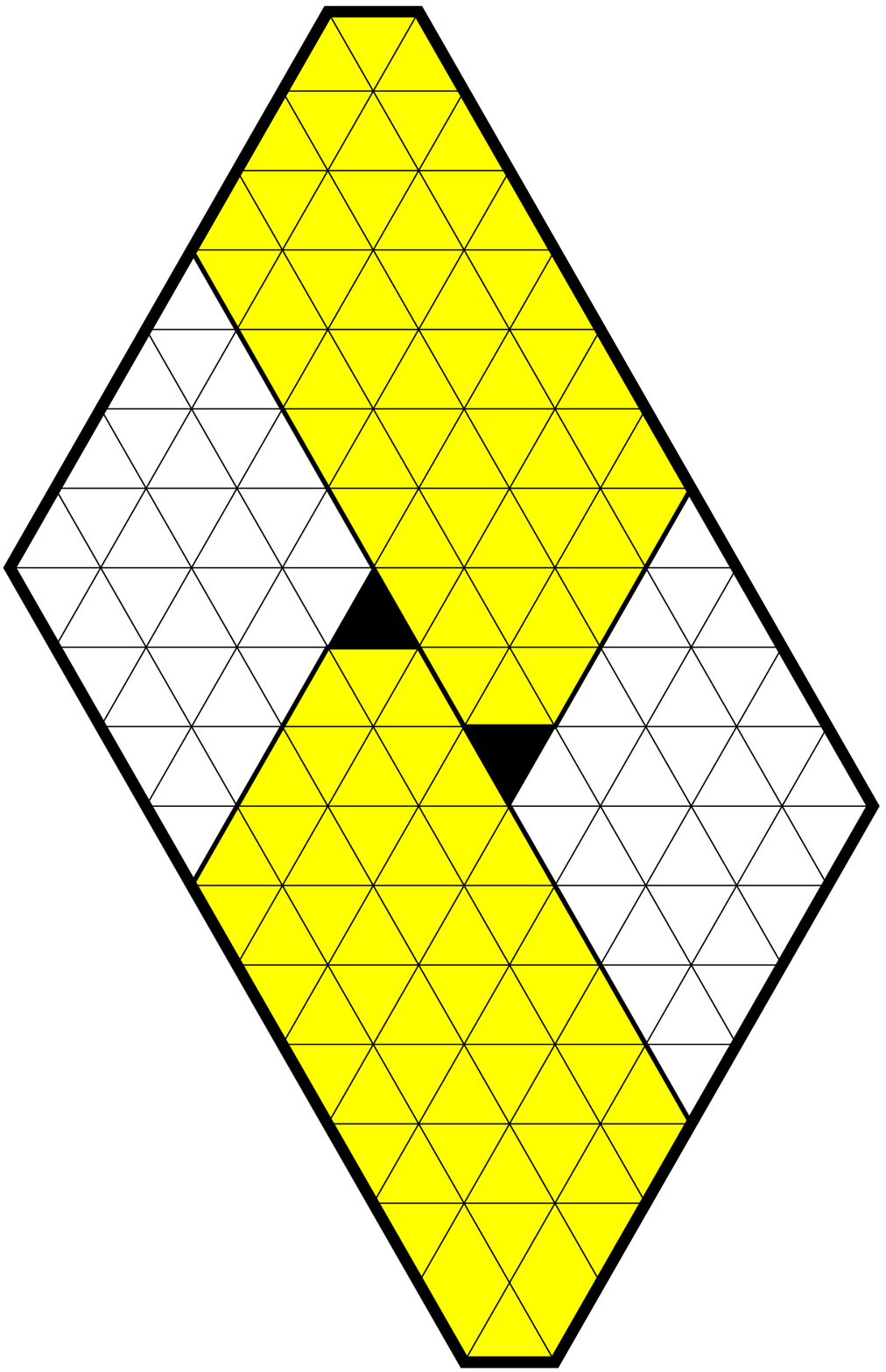}}
\hfill
}
  \caption{\label{fed} The base cases $x=1, k=0$ (left) and $x=1, k=1$ (right).}
\end{figure}

Recall that lozenge tilings of regions on the triangular lattice are in bijection with families of non-intersecting paths of rhombi connecting lattice segments on the boundary of a given direction. By this encoding, the lozenge tilings of $H_{1,y,z}$ can be identified with single paths of rhombi connecting the edges of length 1. A lozenge tiling is centrally symmetric if and only if the corresponding path of rhombi is so. Clearly, any such centrally symmetric path must contain the lozenge $L$ whose short diagonal contains the center of $H_{1,y,z}$ (see the picture on the left in Figure \ref{fed}). Furthermore, such a centrally symmetric path of lozenges is uniquely determined by its top half, which can be arbitrary. Since paths of rhombi connecting the top side of the hexagon with the top side of $L$ are in bijection with tilings of $H_{1,y/2,(z-1)/2}$, we obtain that
\begin{equation}
\M_\odot(B'_{1,y,z,0})=\M\left(H_{1,\frac{y}{2},\frac{z-1}{2}}\right).
\end{equation}  
Formula \eqref{ebd} follows then from MacMahon's theorem \eqref{eaa}.

Assume now that $x=1$ and $k=1$. Then $x$ and $k$ have the same parity, and we are in the situation of part (a) of the theorem. Note that the shaded hexagons in the picture on the right in Figure \ref{fed} must be internally tiled. Since in the remaining portions all tiles are uniquely determined, it follows that
\begin{equation}
\M_\odot(B'_{1,y,z,1})=\M\left(H_{1,\frac{y+1}{2},\frac{z-1}{2}}\right);
\end{equation}  

\parindent0pt
using \eqref{eaa}, this is readily seen to agree with the $x=k=1$ specialization of formula \eqref{ebc}. 
\parindent15pt

\begin{figure}[h]
  \centerline{
\hfill
{\includegraphics[width=0.36\textwidth]{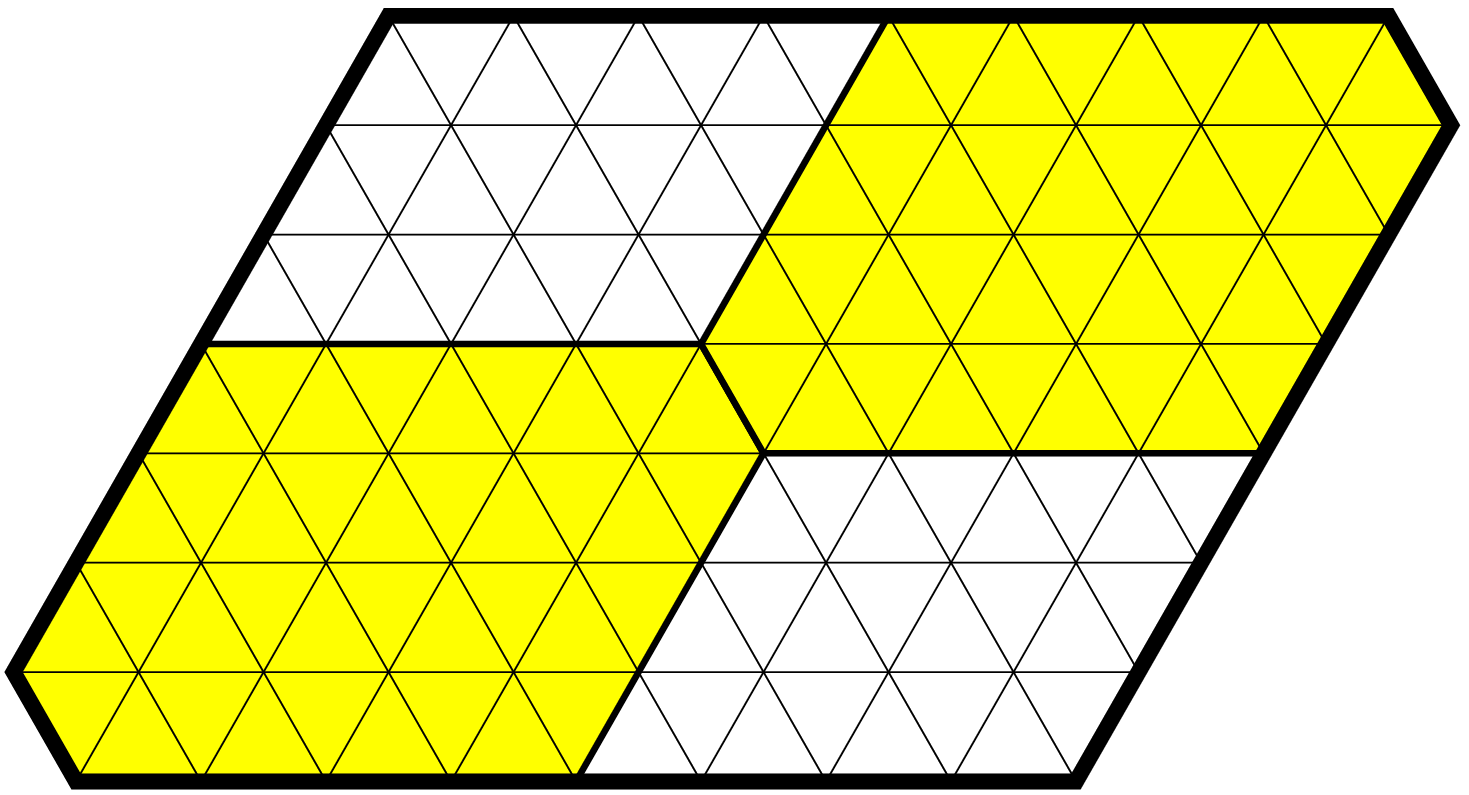}}
\hfill
{\includegraphics[width=0.36\textwidth]{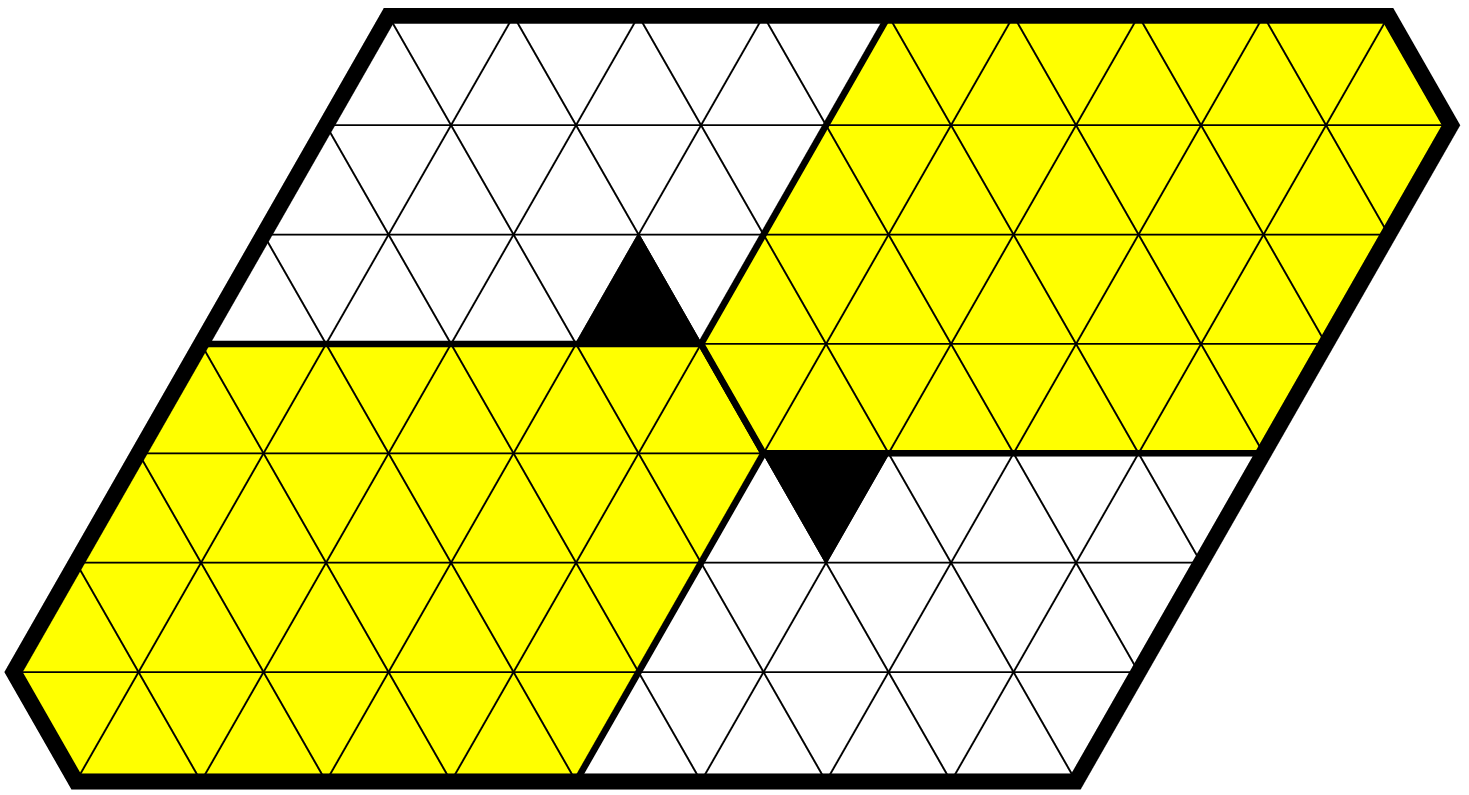}}
\hfill
}
  \caption{\label{fee} The base cases $y=1, k=0$ (left) and $y=1, k=1$ (right).}
\end{figure}

If $y=1$, then either $k=0$ or $k=1$. If $k=0$, then $k$ and $y$ have opposite parity, so we are in the situation of part (a) of the theorem. The region $B'_{x,1,z,0}$ is then just the hexagon $H_{x,1,z}$. One readily sees that the lattice path encoding any centrally symmetric lozenge tiling of $H_{x,1,z}$ must contain the lattice segment containing the center. Then the shaded hexagons on the left in Figure \ref{fee} must be internally tiled. It follows that
\begin{equation}
\M_\odot(B'_{x,1,z,0})=\M\left(H_{\frac{x}{2},1,\frac{z}{2}}\right).
\end{equation}  
Formula \eqref{ebc} follows then from \eqref{eaa}.

Suppose now that $y=1$, $k=1$. Then $y$ and $k$ have the same parity, so we are in the situation covered by part (b) of the theorem. If the path of lozenges $P$ connecting the southwestern side to the northeastern side of $B'_{x,1,z,1}$ does not pass between the two lobes of the removed disconnected bowtie, then the corresponding tiling cannot be centrally symmetric (indeed, the image $P'$ of $P$ through the center is then different from $P$, and has the same starting and ending segments). Therefore centrally symmetric tilings of $B'_{x,1,z,1}$ correspond to paths $P$ that pass through the lobes and are centrally symmetric. Such paths are determined by their portion between the southwestern side and the central segment, which in turn can be identified with the tilings of the lower shaded hexagon in the picture on the right in Figure \ref{fee}. We obtain that
\begin{equation}
\M_\odot(B'_{x,1,z,1})=\M\left(H_{\frac{x}{2},1,\frac{z}{2}}\right),
\end{equation}  
and formula \eqref{ebd} follows from equation \eqref{eaa}.

For part a), there are the additional base cases $k=x$, $k=z$ and $k=y-1$. The first two follow the same way as the analogous base cases in the proof of Theorem \ref{tbb}. In fact, the third also follows the same way, once one observes that in the situation of part a), the unit segment joining the closest vertices in the two removed equilateral triangles is the side of a tile in all tilings of $B'_{x,y,z,k}$. The corresponding base cases for part b) also follow the same way, using this time the observation that, under the assumptions of part b), the lozenge whose short diagonal is the unit segment joining the closest vertices of the two removed equilateral triangles is present in all tilings of $B'_{x,y,z,k}$.

For the induction step, let $x,y,z\geq2$, and assume that formulas \eqref{ebc} and \eqref{ebd} hold for all $B'$-regions for which the sum of their $x$-, $y$- and $z$-parameters are strictly less than $x+y+z$. We need to deduce from this that these formulas also hold for the region $B'_{x,y,z,k}$.

Note that if the indices of $B'_{x,y,z,k}$ are such that $x$, $z$ and $k$ have parity opposite to the parity of $y$, then the same is true for all other regions involved in \eqref{eed}. An analogous statement holds for the case when $x$ and $z$ have parity opposite to the parity of $y$ and $k$. It follows that each of \eqref{ebc} and \eqref{ebd} can be proved separately, by induction on $x+y+z$. The argument is precisely the same as in the proof of Theorem \ref{tbb}. \epf




\section{Concluding remarks}

In this paper we considered centrally symmetric hexagons on the triangular lattice with a bowtie shaped hole (or a shape obtained from a bowtie by separating its lobes one unit) from its center, and we provided simple product formulas that enumerate their centrally symmetric lozenge tilings. This constitutes a generalization of the enumeration of self-complementary plane partitions that fit in a box, which was first proved by Stanley \cite{Sta}.

The proof is based on a new Kuo-style graphical condensation relation that we also prove in this paper. Our extended result is general enough so that this new Kuo-style relation readily affords the proof.

With this, five of the six non-trivial symmetry classes of shamrocks are proved (two were proved in \cite{symffa}, one in \cite{symffb}, and one in \cite{symffc}). The remaining case will be treated in a separate paper.

\end{document}